\documentclass[12pt,reqno]{amsproc}

\usepackage{geometry}
\usepackage{amsthm}

\usepackage{hhline}

\usepackage{mathtools}

\usepackage{mathdots}

\usepackage{color}

\usepackage{pdfsync}

\usepackage{enumitem}

\usepackage{wasysym}

\usepackage{amssymb}

\usepackage{mathrsfs}

\usepackage{bbm}

\DeclareMathAlphabet{\mathpzc}{OT1}{pzc}{m}{it}

\usepackage[all]{xy}

\usepackage{tikz}
\usetikzlibrary{arrows,matrix,decorations.pathmorphing,decorations.pathreplacing,positioning,shapes.geometric,shapes.misc,decorations.markings,decorations.fractals,calc,patterns}

\usepackage{graphicx}

\usepackage{float}

\usepackage[bottom]{footmisc}

\usepackage{moreenum}

\usepackage{calc}

\usepackage{makecell}

\entrymodifiers={+!!<0pt,\fontdimen22\textfont2>}

\usepackage{scalerel,stackengine}
\stackMath
\newcommand\newcheck[1]{%
\savestack{\tmpbox}{\stretchto{%
  \scaleto{%
    \scalerel*[\widthof{\ensuremath{#1}}]{\kern-.6pt\bigwedge\kern-.6pt}%
    {\rule[-\textheight/2]{1ex}{\textheight}}
  }{\textheight}%
}{0.5ex}}%
\stackon[1pt]{#1}{\scalebox{-1}{\tmpbox}}%
}
\stackMath
\newcommand\newhat[1]{%
\savestack{\tmpbox}{\stretchto{%
  \scaleto{%
    \scalerel*[\widthof{\ensuremath{#1}}]{\kern-.6pt\bigwedge\kern-.6pt}%
    {\rule[-\textheight/2]{1ex}{\textheight}}
  }{\textheight}%
}{0.5ex}}%
\stackon[1pt]{#1}{\scalebox{1}{\tmpbox}}%
}

\def\cA{\mathscr{A}}

\def\cC{\mathscr{C}}
\def\cD{\mathscr{D}}
\def\cE{\mathscr{E}}
\def\cF{\mathscr{F}}

\def\cP{\mathscr{P}}

\def\cS{\mathscr{S}}

\def\cW{\mathscr{W}}
\def\cX{\mathscr{X}}
\def\cY{\mathscr{Y}}

\def\BH{\mathbb{H}}

\def\BZ{\mathbb{Z}}

\def\Bk{\mathbbm{k}}

\def\fa{\mathfrak{a}}

\def\fr{\mathfrak{r}}

\mathchardef\mhyphen="2D

\def\adots{\mathinner{\mkern1mu\raise1.0pt\vbox{\kern7.0pt\hbox{.}}\mkern2mu\raise5.0pt\hbox{.}\mkern2mu\raise9.0pt\hbox{.}\mkern1mu}}

\def\ast{{\textstyle *}}

\def\c{\operatorname{c}}

\def\cf{\operatorname{cf}}
\def\Ch{\operatorname{Ch}}

\def\ChN{\operatorname{Ch}_N}

\def\cof{\operatorname{cof}}

\def\dddots{\mathinner{\mkern1mu\raise10.0pt\vbox{\kern7.0pt\hbox{.}}\mkern2mu\raise5.3pt\hbox{.}\mkern2mu\raise1.0pt\hbox{.}\mkern1mu}}
\def\dddotssmall{\mathinner{\mkern1mu\raise7.0pt\vbox{\kern7.0pt\hbox{.}}\mkern-1mu\raise4pt\hbox{.}\mkern-1mu\raise1.0pt\hbox{.}\mkern1mu}}

\def\dual{\operatorname{D}}

\def\Ext{\operatorname{Ext}}

\def\fib{\operatorname{fib}}

\def\Ho{\operatorname{Ho}}
\def\Hom{\operatorname{Hom}}
\def\id{\operatorname{id}}
\def\Image{\operatorname{Im}}

\def\ind{\operatorname{ind}}

\def\Inj{\operatorname{Inj}}

\def\isub{i}
\def\K{\operatorname{K}}
\def\K0{\operatorname{K}_0}

\def\Ker{\operatorname{Ker}}

\def\mod{\operatorname{mod}}
\def\Mod{\operatorname{Mod}}

\def\opp{\operatorname{op}}

\def\perf{\operatorname{perf}}
\def\prj{\operatorname{prj}}
\def\Prj{\operatorname{Prj}}

\def\PSL2{\operatorname{PSL}_2}
\def\psub{p}
\def\Qcpx{Q^{\operatorname{cpx}}}

\def\QNcpx{Q^{N\!\operatorname{-cpx}}}

\def\SL2{\operatorname{SL}_2}

\def\sperf{\operatorname{s.perf}}

\newcommand\Tensor[1]{{\underset{#1}{\otimes}}}

\def\Tor{\operatorname{Tor}}

\def\weq{\operatorname{weq}}

\numberwithin{equation}{section}

\newcommand{\numberseries}{\mdseries}   

\newlength{\thmtopspace}                
\newlength{\thmbotspace}                
\newlength{\thmheadspace}               
\newlength{\thmindent}                  

\newtheoremstyle{bfupright head,slanted body}
                {\thmtopspace}{\thmbotspace}
                {\slshape}{\thmindent}{\bfseries}{.}{\thmheadspace}
                {{\numberseries \thmnumber{{\bf #2} }}\thmnote{#3}}

\newtheoremstyle{bfupright head,upright body}
                {\thmtopspace}{\thmbotspace}
                {\upshape}{\thmindent}{\bfseries}{.}{\thmheadspace}
                {{\numberseries \thmnumber{{\bf #2} }}\thmnote{#3}}

\newtheoremstyle{bfit head,upright body}
                {\thmtopspace}{\thmbotspace}
                {\upshape}{\thmindent}{\upshape}{.}{\thmheadspace}
                {{\numberseries\thmnumber{{\bf #2} }}
                {\bfseries\itshape\thmnote{\negthickspace#3}}}

\newtheoremstyle{it head,upright body}
                {\thmtopspace}{\thmbotspace}
                {\upshape}{\thmindent}{\upshape}{.}{\thmheadspace}
                {{\numberseries\thmnumber{{\bf #2} }}
                {\itshape\thmnote{\negthickspace#3}}}

\newtheoremstyle{fixed bf head,slanted body}
                {\thmtopspace}{\thmbotspace}{\slshape}
                {\thmindent}{\bfseries}{.}{\thmheadspace}
                {{\numberseries \thmnumber{{\bf #2} }}\thmname{#1}\thmnote{ (#3)}}

\newtheoremstyle{fixed bf head,upright body}
                {\thmtopspace}{\thmbotspace}{\upshape}
                {\thmindent}{\bfseries}{.}{\thmheadspace}
                {{\numberseries \thmnumber{{\bf #2} }}\thmname{#1}\thmnote{ (#3)}}

\newtheoremstyle{independent paragraph}
                {\thmtopspace}{\thmbotspace}
                {\upshape}{\thmindent}{\upshape}{}{0pt}
                {\thmnote{#3 }}

\newtheoremstyle{subparagraph}
                {\thmbotspace}{\thmbotspace}
                {\upshape}{\thmindent}{\upshape}{}{0pt}
                {\thmnote{#3 }}

\newtheoremstyle{notes}
                {\thmtopspace}{\thmbotspace}
                {\ttfamily}{\thmindent}{\ttfamily\small }{}{0pt}
                {\thmnote{#3 }}
                
\newtheoremstyle{numbered paragraph}
                {\thmtopspace}{\thmbotspace}{\upshape}
                {\thmindent}{\upshape}{}{\thmheadspace}
                {{\numberseries \thmnumber{\bf #2.}}}

\theoremstyle{bfupright head,slanted body}
\newtheorem{res}{}[section]             \newtheorem*{res*}{}

\theoremstyle{bfit head,upright body}
                 \newtheorem*{com*}{}

\theoremstyle{bfupright head,upright body}
\newtheorem{bfhpg}[res]{}               \newtheorem*{bfhpg*}{}

\theoremstyle{it head,upright body}
               \newtheorem*{ithpg*}{}

\theoremstyle{fixed bf head,slanted body}
\newtheorem{Theorem}[res]{Theorem}          \newtheorem*{Theorem*}{Theorem}
      \newtheorem*{Proposition*}{Proposition}
        \newtheorem*{Corollary*}{Corollary}
            \newtheorem*{Lemma*}{Lemma}

\theoremstyle{fixed bf head,upright body}
       \newtheorem*{Definition*}{Definition}
      \newtheorem*{Observation*}{Observation}
\newtheorem{Remark}[res]{Remark}           \newtheorem*{Remark*}{Remark}
          \newtheorem*{Example*}{Example}
         \newtheorem*{Exercise*}{Exercise}
\newtheorem{Setup}[res]{Setup}            \newtheorem{Setup*}{Setup}
         \newtheorem{Notation*}{Notation}
     \newtheorem{Construction*}{Construction}
\newtheorem{Conjecture}[res]{Conjecture}     \newtheorem{Conjecture*}{Conjecture}

\theoremstyle{numbered paragraph}

\theoremstyle{subparagraph}

\theoremstyle{notes}

  {\begin{list}{}{%
    \settowidth{\labelwidth}{\textbf{#1:}}%
    \setlength{\leftmargin}{\labelwidth}\addtolength{\leftmargin}{\labelsep}}}%
  {\end{list}}

\makeatletter
\@namedef{subjclassname@2020}{%
  \textup{2020} Mathematics Subject Classification}
\makeatother

\begin{document}

\setlength{\parindent}{0pt}
\setlength{\parskip}{7pt}

\title[The $Q$-shaped derived category]{A brief introduction to the $Q$-shaped derived category}

\author{Henrik Holm}

\address{Department of Mathematical Sciences, Universitetsparken 5, University of Copenhagen, 2100 Copenhagen {\O}, Denmark} 
\email{holm@math.ku.dk}

\urladdr{http://www.math.ku.dk/\~{}holm/}

\author{Peter J\o rgensen}

\address{Department of Mathematics, Aarhus University, Ny Munkegade 118, 8000 Aarhus C, Denmark}
\email{peter.jorgensen@math.au.dk}

\urladdr{https://sites.google.com/view/peterjorgensen}





\begin{abstract} 

A chain complex can be viewed as a representation of a certain quiver with relations, $\Qcpx$.  The vertices are the integers, there is an arrow $q \xrightarrow{} q-1$ for each integer $q$, and the relations are that consecutive arrows compose to $0$.  Hence the classic derived category $\cD$ can be viewed as a category of representations of $\Qcpx$.  

\bigskip
\noindent
It is an insight of Iyama and Minamoto that the reason $\cD$ is well behaved is that, viewed as a small category, $\Qcpx$ has a Serre functor.  Generalising the construction of $\cD$ to other quivers with relations which have a Serre functor results in the {\em $Q$-shaped derived category}, $\cD_Q$. 

\bigskip
\noindent
Drawing on methods of Hovey and Gillespie, we developed the theory of $\cD_Q$ in three recent papers.  This paper offers a brief introduction to $\cD_Q$, aimed at the reader already familiar with the classic derived category.

\end{abstract}

\maketitle

\tableofcontents

\setcounter{section}{-1}
\section{Introduction}
\label{sec:introduction}

The classic derived category $\cD( A )$ of a ring provides a framework for homological algebra.  Its objects are chain complexes, which can also be viewed as representations of a certain quiver with relations, $\Qcpx$, defined by Figure \ref{fig:linear_quiver} with the relations that consecutive arrows compose to $0$.

It is an insight of Iyama and Minamoto that the key property making $\cD( A )$ work is that, viewed as a small category, $\Qcpx$ has a Serre functor; see \cite{Iyama-Minamoto-1}, \cite[sec.\ 2]{Iyama-Minamoto-2}.  Generalising the construction to other quivers with relations which have a Serre functor, or more generally to any suitable category $Q$ which has a Serre functor, results in the {\em $Q$-shaped derived category} $\cD_Q( A )$.  See Setup \ref{set:blanket} for the precise conditions imposed on $Q$.

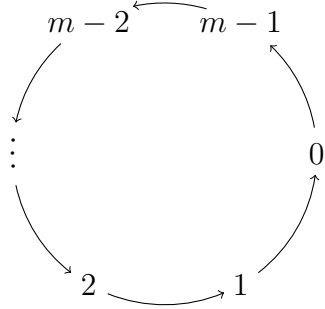
\begin{figure}
\begin{tikzpicture}[scale=2]
  \node at (0:1.0){$0$};
  \draw[->] (10:1.0) arc (10:46:1.0);
  \node at (60:1.0){$m-1$};
  \draw[->] (74:1.0) arc (74:102:1.0);
  \node at (120:1.0){$m-2$};
  \draw[->] (133:1.0) arc (133:168:1.0);
  \node at (175:1.0){$\cdot$};
  \node at (180:1.0){$\cdot$};
  \node at (185:1.0){$\cdot$};  
  \draw[->] (192:1.0) arc (192:232:1.0);
  \node at (240:1.0){$2$};
  \draw[->] (248:1.0) arc (248:293:1.0);
  \node at (300:1.0){$1$};
  \draw[->] (308:1.0) arc (308:352:1.0);
\end{tikzpicture}
\caption{A cyclic quiver with $m$ vertices.}
\label{fig:cyclic_quiver}
\end{figure}

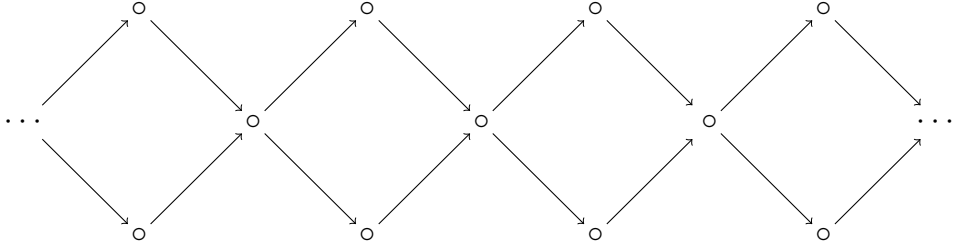
\begin{figure}
\begin{tikzpicture}[scale=1.5]
  \node at (-3,0){$\cdots$};
  \draw[->] (-2.85,-0.15) to (-2.1,-0.9);
  \draw[->] (-2.85,0.15) to (-2.1,0.9);  
  \node at (-2,-1){$\circ$};
  \node at (-2,1){$\circ$};
  \draw[->] (-1.9,-0.9) to (-1.1,-0.1);  
  \draw[->] (-1.9,0.9) to (-1.1,0.1);  
  \node at (-1,0){$\circ$};
  \draw[->] (-0.9,-0.1) to (-0.1,-0.9);
  \draw[->] (-0.9,0.1) to (-0.1,0.9);
  \node at (0,-1){$\circ$};
  \node at (0,1){$\circ$};  
  \draw[->] (0.1,-0.9) to (0.9,-0.1);    
  \draw[->] (0.1,0.9) to (0.9,0.1);    
  \node at (1,0){$\circ$};
  \draw[->] (1.1,-0.1) to (1.9,-0.9);    
  \draw[->] (1.1,0.1) to (1.9,0.9);    
  \node at (2,-1){$\circ$};
  \node at (2,1){$\circ$};
  \draw[->] (2.1,-0.9) to (2.85,-0.15);    
  \draw[->] (2.1,0.9) to (2.85,0.15);    
  \node at (3,0){$\circ$};
  \draw[->] (3.1,-0.1) to (3.9,-0.9);    
  \draw[->] (3.1,0.1) to (3.9,0.9);    
  \node at (4,-1){$\circ$};
  \node at (4,1){$\circ$};
  \draw[->] (4.1,-0.9) to (4.85,-0.15);    
  \draw[->] (4.1,0.9) to (4.85,0.15);    
  \node at (5,0){$\cdots$};    
\end{tikzpicture}
\caption{The repetitive quiver $\BZ A_3$.}
\label{fig:ZA3}
\end{figure}

The $Q$-shaped derived category shares several attractive properties of $\cD( A )$; for instance, it is a compactly generated triangulated category.  At the same time, varying $Q$ provides the freedom to construct bespoke triangulated categories.  For instance, if $Q = \QNcpx$ is defined by Figure \ref{fig:linear_quiver} with the relations that any $N$ consecutive arrows compose to $0$, then $\cD_Q( A )$ is the derived category of $N$-complexes introduced in \cite{Iyama-Kato-Miyachi-N}.  If $Q$ is defined by Figure \ref{fig:cyclic_quiver} with the relations that consecutive arrows compose to $0$, then $\cD_Q( A )$ is the derived category of $m$-periodic complexes, which has the special feature that $\Sigma^{ 2m } \cong \id$ where $\Sigma$ is the suspension functor.  (The power $2m$ can be replaced by $m$ if $m$ is even, but if $m$ is odd, then $\Sigma^m$ flips the sign of the differential.)  One can also pick a more complicated $Q$, e.g.\ defined by Figure \ref{fig:ZA3} with mesh relations.

The construction of $\cD_Q( A )$ draws heavily on methods of Hovey and Gillespie, in particular \cite{Hovey-MathZ} and \cite{Gillespie-BLMS}.  See \cite{Gillespie-construct}, \cite{Gillespie-exact}, \cite{Gillespie-flat}, \cite{Hovey-book} for additional background.

This paper offers a brief introduction to the theory of $\cD_Q( A )$, which was developed in the papers \cite{HJ-Model_cats}, \cite{HJ-JLMS}, \cite{HJ-arXiv}.  Most sections begin with one or more items of abstract theory, followed by concrete implementations in the following standing examples.
\begin{itemize}
\setlength\itemsep{4pt}

  \item  The category $Q = \Qcpx$ defined by Figure \ref{fig:linear_quiver} with the relations that consecutive arrows compose to $0$.  Here $\cD_Q( A ) = \cD( A )$ is the classic derived category.  (Up to now, $\Qcpx$ denoted a quiver with relations.  It will henceforth denote the corresponding small category; see \ref{bfhpg:examples_of_Q}.)
  
  \item  The category $Q = \QNcpx$ defined by Figure \ref{fig:linear_quiver} with the relations that any $N$ consecutive arrows compose to $0$ for a fixed integer $N \geqslant 2$.  Here $\cD_Q( A ) = \cD_N( A )$ is the derived category of $N$-complexes.

\end{itemize}
We end the introduction with a preview, which on its own can serve as an even briefer introduction to $\cD_Q( A )$.  The paper is divided into Parts \ref{part:Frobenius} through \ref{part:compendium}, and the preview covers Parts \ref{part:Frobenius} through \ref{part:compact}, comprising Sections \ref{sec:H_E_weq} through \ref{sec:compact}.  Part \ref{part:compendium} contains two appendices on some key classes of categories: Frobenius, triangulated, and abelian model categories.

Let	$\Bk$ be a hereditary noetherian commutative ring, $A$ a $\Bk$-algebra, and $Q$ a category satisfying the conditions in Setup \ref{set:blanket}, with $Q_0$ denoting the class of objects of $Q$.  Note that $Q$ is often defined by a quiver with relations; see \ref{bfhpg:examples_of_Q}.  Let ${}_{ Q,A }\!\Mod$ be the abelian category of $\Bk$-linear functors $Q \xrightarrow{} \Mod( A )$ where $\Mod( A )$ is the category of $A$-left modules.

\begin{bfhpg}[The Frobenius approach to $\cD_Q( A )$ (preview of Part \ref{part:Frobenius})]
\label{bfhpg:preview_A}
This part constructs $\cD_Q( A )$ in two different ways as the stable category of a Frobenius category.

{\em Section \ref{sec:H_E_weq}:} For each $q \in Q_0$ and integer $i \geqslant 0$, there are {\em (co)homology functors}
\[  
  \BH_{ [q] }^i \:,\: \BH^{ [q] }_i:
  {}_{ Q,A }\!\Mod \xrightarrow{} \Mod( A ).
\]
  For $i = 1$ they generalise the classic homology functors $H_j: \Ch( A ) \xrightarrow{} \Mod( A )$ where $\Ch( A )$ is the category of chain complexes over $\Mod( A )$.  The full subcategory of {\em exact} objects in ${}_{ Q,A }\!\Mod$ is
\[
  \cE = \{ X \in {}_{ Q,A }\!\Mod \mid \mbox{$\BH^1_{ [q] }( X ) = 0$ for each $q \in Q_0$} \};
\]
it generalises the exact complexes.  The {\em weak equivalences} in ${}_{ Q,A }\!\Mod$ are
\[
  \weq = 
  \Bigg\{ \mbox{morphisms $\varphi$ in ${}_{ Q,A }\!\Mod$}
  \:\Bigg|\! 
  \begin{array}{l}
    \mbox{$\BH^1_{ [q] }( \varphi )$ and $\BH^2_{ [q] }( \varphi )$ are} \\[1mm]
    \mbox{isomorphisms for each $q \in Q_0$} 
  \end{array}
  \Bigg\};
\]
they generalise the quasiisomorphisms of complexes.  It may be unexpected that $\BH^1_{ [q] }( \varphi )$ and $\BH^2_{ [q] }( \varphi )$ are both required to be isomorphisms in the formula, but this is necessary for the theory to work.

{\em Section \ref{sec:E_perp}:}  Let ${}_{ Q,A }\!\Prj$ and ${}_{ Q,A }\!\Inj$ be the full subcategories of projective, respectively injective objects of ${}_{ Q,A }\!\Mod$.  The {\em semiprojective} objects of ${}_{ Q,A }\!\Mod$ are the objects of ${}^{ \perp }\cE$, which is a Frobenius category with ${}_{ Q,A }\Prj$ as its class of projective-injective objects.  Here $\perp$ indicates a perpendicular full subcategory with respect to $\Ext^1$, see \ref{bfhpg:QAMod}.  Semiprojective objects generalise semiprojective complexes.  Similarly, the {\em semiinjective} objects of ${}_{ Q,A }\!\Mod$ are the objects of $\cE^{ \perp }$, which is a Frobenius category with ${}_{ Q,A }\Inj$ as its class of projective-injective objects.  Semiinjective objects generalise semiinjective complexes.

{\em Section \ref{sec:DQ}:} The {\em $Q$-shaped derived category of $A$} is obtained from ${}_{ Q,A }\!\Mod$ by formally inverting each weak equivalence,
\[
  \cD_Q( A ) = \weq^{ -1 }\!{}_{ Q,A }\!\Mod.
\]
There are equivalences of categories
\[
\tag{$\ast$}
  \frac{ {}^{ \perp }\cE }{ {}_{ Q,A }\Prj }
  \cong \cD_Q( A ) \cong
  \frac{ \cE^{ \perp } }{ {}_{ Q,A }\Inj }.
\]
Here $\frac{ {}^{ \perp }\cE }{ {}_{ Q,A }\Prj }$ and $\frac{ \cE^{ \perp } }{ {}_{ Q,A }\Inj }$ are the stable categories of the Frobenius categories ${}^{ \perp }\cE$ and $\cE^{ \perp }$.  Hence they are triangulated categories, and in Part \ref{part:Frobenius} of the paper, we view them as the de facto definition of $\cD_Q( A )$.  The equivalences permit the concrete computation of $\Hom$ spaces in $\cD_Q( A )$.  For instance,
\[
  \Hom_{ \cD_Q( A ) }( X,X' ) \cong 
  \frac{ \Hom_{ {}_{ Q,A }\!\Mod }( P,P' ) }{\{\mbox{ morphisms which factorise through a projective object }\}}.
\]
Here $P$ and $P'$ are semiprojective resolutions of $X$ and $X'$; that is, $P$ and $P'$ are semiprojective objects with weak equivalences $P \xrightarrow{} X$ and $P' \xrightarrow{} X'$.

{\em Section \ref{sec:suspension}:} Gives sample computations of the suspension functor of $\cD_Q( A )$.
\end{bfhpg}

\begin{bfhpg}[The model category approach to $\cD_Q( A )$ (preview of Part \ref{part:model})]
This part constructs the projective and injective model category structures on ${}_{ Q,A }\!\Mod$ and obtains $\cD_Q( A )$ as the corresponding homotopy category where each weak equivalence has been formally inverted.

{\em Section \ref{sec:cotorsion}:} There are cotorsion pairs
\[  
  ( {}^{ \perp }\cE,\cE ) \;,\;
  ( {}_{ Q,A }\!\Prj,{}_{ Q,A }\!\Mod) \;,\;
  ( \cE,\cE^{ \perp } ) \;,\; 
  ( {}_{ Q,A }\!\Mod,{}_{ Q,A }\!\Inj )
\]
in ${}_{ Q,A }\!\Mod$.
  
{\em Section \ref{sec:model}:} There are two model category structures on ${}_{ Q,A }\!\Mod$: The {\em projective} and the {\em injective model category structures}.  They arise by applying Hovey's Theorem (Theorem \ref{thm:Hovey} in Appendix \ref{app:II}) to the so-called Hovey triples $( {}^{ \perp }\cE,\cE,{}_{ Q,A }\!\Mod )$ and $( {}_{ Q,A }\!\Mod,\cE,\cE^{ \perp } )$, which are obtained from the cotorsion pairs of Section \ref{sec:cotorsion}.  Both model category structures have the class $\weq$ from \ref{bfhpg:preview_A} as their weak equivalences.  Hence they have the same homotopy category:
\[
  \Ho( {}_{ Q,A }\!\Mod )
  = \weq^{ -1 }\!{}_{ Q,A }\!\Mod
  = \cD_Q( A ).
\]
The equivalences $(\ast)$ from \ref{bfhpg:preview_A} are now obtained from a theorem by Gillespie (Theorem \ref{thm:Gillespie}).
\end{bfhpg}

\begin{bfhpg}[Compact, perfect, and strictly perfect objects in $\cD_Q( A )$ (preview of Part \ref{part:compact})]
This part presents some classes of objects and some properties of $\cD_Q( A )$. 

{\em Section \ref{sec:compact}:} $\cD_Q( A )$ has full subcategories $\cD_Q^{ \c }( A )$, $\cD_Q^{ \perf }( A )$, $\cD_Q^{ \sperf }( A )$ of {\em compact}, {\em perfect}, and {\em strictly perfect} objects, which enjoy certain relations.  A key property is that $\cD_Q( A )$ is a compactly generated triangulated category generated by ``stalk functors'', which send one object of $Q_0$ to $A$ and all other objects to $0$. 
\end{bfhpg}

\begin{figure}
\begin{tikzpicture}[scale=1.5]
  \node at (-3,0){$\cdots$};
  \draw[->] (-2.75,0) to (-2.2,0);
  \node at (-2,0){$2$};
  \draw[->] (-1.8,0) to (-1.2,0);  
  \node at (-1,0){$1$};
  \draw[->] (-0.8,0) to (-0.2,0);    
  \node at (0,0){$0$};
  \draw[->] (0.2,0) to (0.75,0);    
  \node at (1,0){$-1$};
  \draw[->] (1.25,0) to (1.75,0);    
  \node at (2,0){$-2$};    
  \draw[->] (2.25,0) to (2.75,0);    
  \node at (3,0){$\cdots$};    
\end{tikzpicture}
\caption{The linear quiver which underlies chain complexes and $N$-complexes.}
\label{fig:linear_quiver}
\end{figure}

\section{Preliminaries}

\begin{Setup}
\label{set:blanket}
This paragraph states the setup which will be assumed in the rest of the paper.  Many of the results require less than the full setup, but we refer the reader to the original papers for specifics.
\begin{itemize}
\setlength\itemsep{4pt}

  \item  $\Bk$ is a hereditary noetherian commutative ring.  
  
  \item  $A$ is a $\Bk$-algebra.  
  
  \item  $\Mod( \Bk )$ is the category of $\Bk$-modules, and $\Mod( A )$ is the category of $A$-left modules.

\end{itemize}
Typical examples of $\Bk$ are $\BZ$ or a field.  If $\Bk = \BZ$, then $A$ is just a ring.  Moreover, 
\begin{itemize}
\setlength\itemsep{4pt}

  \item  $Q$ is a category satisfying the following conditions, mainly due to \cite[thm.\ 1.6]{DSS}.  Here $Q_0$ denotes the class of objects of $Q$, and $Q( -,- )$ denotes the homomorphism functor of $Q$.

\end{itemize}

{\em Preadditivity:}  $Q$ is $\Bk$-preadditive; that is, each $\Hom$ set is equipped with a structure of $\Bk$-module, and composition of morphisms is $\Bk$-bilinear.

{\em Hom finiteness:}  Each $\Hom$ set in $Q$ is a finitely generated free $\Bk$-module.
  
{\em Local boundedness:}  For each $q \in Q_0$, the following sets are finite.
\[
  \{ p \in Q_0 \mid Q( q,p ) \neq 0 \}
  \;\;,\;\;
  \{ p \in Q_0 \mid Q( p,q ) \neq 0 \}
\]

{\em Serre functor:}  There exists a $\Bk$-linear automorphism $S$ of $Q$, called the Serre functor, such that there are isomorphisms
\[
  Q( p,q ) \cong \Hom_{ \Bk }\!\big( Q( q,Sp ),\Bk \big)
\]
which are natural in $p,q \in Q_0$.

{\em Strong retraction:}  $Q$ is equipped with decompositions of $\Bk$-modules
\[  
  Q( q,q ) = (\Bk \cdot \id_q) \oplus \fr_q
\]  
for $q \in Q_0$ such that
\medskip
\begin{enumerate}
\setlength\itemsep{4pt}

  \item  $\fr_q \circ \fr_q \subseteq \fr_q$,
  
  \item  $Q( p,q ) \circ Q( q,p ) \subseteq \fr_q$ for $p \neq q$.

\end{enumerate}

{\em Nilpotence:}  The ideal in $Q$ defined by
\[
  \fr( p,q )
  =
  \left\{
    \begin{array}{cl}
      \fr_q    & \mbox{ if $p=q$, }\\[1mm]
      Q( p,q ) & \mbox{ if $p \neq q$ }
    \end{array}
  \right.
\]
for $p,q \in Q_0$ is called the {\em pseudoradical}.  It must satisfy $\fr^N = 0$ for some integer $N \geqslant 1$.

These conditions are self dual in the sense that they hold for $Q$ and $Q^{ \opp }$ simultaneously.
\end{Setup}

\begin{bfhpg}[Remarks on condition ``Strong retraction'']
\label{bfhpg:Strong_retraction}
\hspace{0.1mm}
\begin{itemize}
\setlength\itemsep{4pt}

  \item  Condition ``Strong retraction'' requires a fixed, global choice of decompositions $Q( q,q ) = (\Bk \cdot \id_q) \oplus \fr_q$, and the pseudoradical $\fr$ depends on the choice.

  \item  In condition ``Strong retraction'', note that ``$p \neq q$'' means ``not equal'' as opposed to ``not isomorphic''.  This implies that different objects of $Q$ are non-isomorphic; see \cite[rmk.\ 7.6]{HJ-JLMS}.
  
\end{itemize}
\end{bfhpg}

\begin{bfhpg}[Cycles in $Q$]
\label{bfhpg:cycle}
A {\em cycle} in the category $Q$ is a diagram in $Q$,
\[
  q_1 \xrightarrow{} q_2 \xrightarrow{} \cdots \xrightarrow{} q_{ n-1 } \xrightarrow{} q_n,
\]
with $q_n = q_1$ where each morphism is non-zero and belongs to the pseudoradical $\fr$.  Note that this definition depends on the fixed, global choice of decompositions $Q( q,q ) = (\Bk \cdot \id_q) \oplus \fr_q$; see \ref{bfhpg:Strong_retraction}.
\end{bfhpg}

\begin{bfhpg}[The categories ${}_{ Q }\!\Mod$ and ${}_{ Q,A }\!\Mod$]
\label{bfhpg:QAMod}
The following categories will be used frequently.
\begin{itemize}
\setlength\itemsep{4pt}

  \item  ${}_{ Q }\!\Mod
  = \{ \mbox{ $\Bk$-linear functors $Q \xrightarrow{} \Mod( \Bk )$ } \}$.  We think of this as the category of $\Mod( \Bk )$-valued representations of $Q$.
  
  \item  ${}_{ Q,A }\!\Mod
  = \{ \mbox{ $\Bk$-linear functors $Q \xrightarrow{} \Mod( A )$ } \}$.  We think of this as the category of $\Mod( A )$-valued representations of $Q$.  

\end{itemize}
We list some properties of ${}_{ Q }\!\Mod$ and ${}_{ Q,A }\!\Mod$.  
\begin{itemize}
\setlength\itemsep{4pt}

  \item  They are abelian; indeed, they are Grothendieck categories by \cite[prop.\ 3.12]{HJ-JLMS}.  
  
  \item  They have enough projective objects by \cite[prop.\ 3.12(a)]{HJ-JLMS} and enough injective objects by \cite[thm.\ 1.10.1]{Grothendieck-Tohoku}.  
    
  \item  They have set indexed limits and colimits by \cite[chp.\ V, preamble and cor.\ X.4.4]{Stenstroem}.

\end{itemize}
Note that \cite[prop.\ 3.12(a)]{HJ-JLMS} even implies that ${}_{ Q }\!\Mod$ and ${}_{ Q,A }\!\Mod$ have projective generators.

The following notation will be used.
\begin{itemize}
\setlength\itemsep{4pt}

  \item  The full subcategories of projective and injective objects of ${}_{ Q,A }\!\Mod$ are denoted by ${}_{ Q,A }\!\Prj$ and ${}_{ Q,A }\!\Inj$.  

  \item  The $\Hom$ and $\Ext$ functors of ${}_{ Q }\!\Mod$ are denoted by $\Hom_{ Q }$ and $\Ext_{ Q }^i$. 

  \item  The $\Hom$ and $\Ext$ functors of ${}_{ Q,A }\!\Mod$ are denoted by $\Hom_{ Q,A }$ and $\Ext_{ Q,A }^i$.  

  \item  If $\cY$ is a class of objects in an abelian category $\cA$, then its left and right perpendicular full subcategories are
\[
  {}^{ \perp }\cY
  = \{ X \in \cA \mid \Ext_{ \cA }^1( X,\cY ) = 0 \} \;\;,\;\;
  \cY^{ \perp }
  = \{ Z \in \cA \mid \Ext_{ \cA }^1( \cY,Z ) = 0 \}.
\]

\end{itemize}
Finally, we also need tensor and $\Tor$ functors.
\begin{itemize}
\setlength\itemsep{4pt}

  \item  There is a tensor product
\[
  - \Tensor{Q} - : 
  {}_{ Q^{ \opp } }\!\Mod \times {}_{ Q }\!\Mod
  \xrightarrow{}
  \Mod( \Bk ),
\]
see \cite[p.\ 93]{Oberst-Roehrl}.  Its $i$th left derived functor is
\[
  \Tor^Q_i( -,- ) : 
  {}_{ Q^{ \opp } }\!\Mod \times {}_{ Q }\!\Mod
  \xrightarrow{}
  \Mod( \Bk ).
\]

\end{itemize}
\end{bfhpg}

\begin{bfhpg}[Examples of $Q$]
\label{bfhpg:examples_of_Q}
The category $Q$ is often defined by a quiver with relations, whose category of $\Mod( A )$-valued representations can then be identified with ${}_{ Q,A }\!\Mod$.  See \cite[secs.\ II.1 and II.2]{ASS} for background on quivers with relations.  Given a quiver, the vertices can be declared the objects of a $\Bk$-preadditive category, and the $\Bk$-linear combinations of paths can be declared the morphisms.  Dividing by the ideal $\fa$ defined by an admissible set of $\Bk$-linear relations gives a new $\Bk$-preadditive category, and this is a candidate for $Q$.  A pseudoradical is defined by setting
\[
  \fr( p,q )
  =
  \{ \mbox{ $\Bk$-linear combinations of paths from $p$ to $q$ of length $\geqslant 1$ } \}/\fa( p,q ).
\]

The conditions in Setup \ref{set:blanket} apart from ``Preadditivity'' are far from automatic.  In particular, the existence of a Serre functor is fairly special.

However, there are many examples where Setup \ref{set:blanket} is satisfied, and we mention the following in particular.
\begin{itemize}
\setlength\itemsep{4pt}

  \item  $Q = \Qcpx$ is defined by Figure \ref{fig:linear_quiver} with the relations that consecutive arrows compose to $0$.  Then ${}_{ Q,A }\!\Mod = \Ch( A )$ where $\Ch( A )$ is the category of chain complexes and chain maps over $\Mod( A )$.  The action of the Serre functor on objects is $S( q ) = q-1$.
  
  \item  $Q = \QNcpx$ is defined by Figure \ref{fig:linear_quiver} with the relations that any $N$ consecutive arrows compose to $0$ for a fixed integer $N \geqslant 2$.  Then ${}_{ Q,A }\!\Mod = \ChN( A )$ where $\ChN( A )$ is the category of $N$-complexes and morphisms of $N$-complexes over $\Mod( A )$; see \cite[secs.\ 0 and 1]{Kapranov}.  The action of the Serre functor on objects is $S( q ) = q-N+1$.

  \item  Let $Q$ be defined by Figure \ref{fig:cyclic_quiver} with the relations that consecutive arrows compose to $0$.  Here $m \geqslant 1$ is a fixed integer.  Then ${}_{ Q,A }\!\Mod$ can be identified with the category of $m$-periodic chain complexes and chain maps over $\Mod( A )$.  The action of the Serre functor on objects is $S( q ) = (q-1) \bmod m$.

  \item  Let $Q$ be defined by Figure \ref{fig:ZA3} modulo mesh relations; see \cite[sec.\ 0.vii]{HJ-Model_cats}.  Then ${}_{ Q,A }\!\Mod$ is a category not mentioned in standard textbooks.  The action of the Serre functor on objects is given by reflection in a central horizontal line through Figure \ref{fig:ZA3} followed by translation by one vertex to the right.

\end{itemize}
Finally, we give an example which is not based on a quiver.  Assume that $\Bk$ is a field and let $\Lambda$ be a finite dimensional self injective $\Bk$-algebra.
\begin{itemize}
\setlength\itemsep{4pt}

  \item  Let $Q$ be a skeleton of $\ind( \prj \Lambda )$, the category of indecomposable finitely generated projective $\Lambda$-left modules.  Then ${}_{ Q,A }\!\Mod$ can be identified with the category of $\Lambda^{ \opp } \Tensor{\Bk} A$-left modules.  The Serre functor is given by $S( - ) = \dual\!\Lambda \Tensor{\Lambda} -$ where $\dual\!\Lambda = \Hom_{ \Bk }( \Lambda,\Bk )$.  Note that the tensor product must be chosen with values in the skeleton $Q$ of $\ind( \prj \Lambda )$.

\end{itemize}
\end{bfhpg}

\part{The Frobenius approach to $\cD_Q( A )$}
\label{part:Frobenius}

This part constructs $\cD_Q( A )$ in two different ways as the stable category of a Frobenius category.

\section{The functors $\BH_{ [q] }^i$ and $\BH^{ [q] }_i$, the class $\cE$, and the class $\weq$}
\label{sec:H_E_weq}

\begin{bfhpg}[The functors $\BH_{ [q] }^i$ and $\BH^{ [q] }_i$]
\label{bfhpg:the_functors_H}
Let $q \in Q_0$ be given.  The {\em stalk functors at $q$} are defined as follows; see \cite[prop.\ 7.15]{HJ-JLMS}.
\begin{itemize}
\setlength\itemsep{4pt}

  \item $S \langle q \rangle = Q( q,- )/\fr( q,- )$
  
  \item $S \{ q \} = Q( -,q )/\fr( -,q )$

\end{itemize}
They are objects of ${}_{ Q }\!\Mod$, respectively ${}_{ Q^{ \opp } }\!\Mod$; that is, they are $\Bk$-linear functors $Q \xrightarrow{} \Mod( \Bk )$, respectively $Q^{ \opp } \xrightarrow{} \Mod( \Bk )$.  If $p \in Q_0$ then
\[
  S \langle q \rangle ( p )
  =
  \left\{
    \begin{array}{cl}
      \Bk & \mbox{ for } q=p, \\[1mm]
      0   & \mbox{ for } q \neq p
    \end{array}
  \right.
\]
by \cite[lem.\ 7.10]{HJ-JLMS}.  The functor $S \langle q \rangle$ generalises the simple representation of $Q$ at $q$ known from quiver representation theory; see \cite[def.\ 2.2(a)]{Schiffler}.  The functor $S \{ q \}$ generalises the simple representation of $Q^{ \opp }$ at $q$.

Now let $q \in Q_0$ and $i \in \BZ$ be given.  Recalling the functors $\Ext_Q^i$ and $\Tor^Q_i$ from \ref{bfhpg:QAMod}, the {\em $i$'th (co)homology functors at $q$} are defined as follows; see \cite[def.\ 7.11]{HJ-JLMS}.
\begin{itemize}
\setlength\itemsep{4pt}
\setcounter{enumi}{4}

  \item  $\BH^i_{ [q] }( - ) = \Ext^i_Q( S \langle q \rangle,- )$
  
  \item  $\BH_i^{ [q] }( - ) = \Tor_i^Q( S \{ q \}, - )$
  
\end{itemize}
They are $\Bk$-linear functors ${}_{ Q,A }\!\Mod \xrightarrow{} \Mod( A )$.  The values have $A$-structures induced by the $A$-structures of the arguments of the functors.
\end{bfhpg}

\begin{bfhpg}[The class $\cE$ of exact objects]
\label{bfhpg:E}
The full subcategory of {\em exact objects} in ${}_{ Q,A }\!\Mod$ is
\begin{itemize}
\setlength\itemsep{4pt}
\setcounter{enumi}{1}

  \item
\vspace{-4ex}
$
  \begin{array}{rcl}
  \phantom{x} && \\[2mm]
  \cE & = & \{ X \in {}_{ Q,A }\!\Mod \mid \mbox{$\BH^1_{ [q] }( X ) = 0$ for each $q \in Q_0$} \} \\[2mm]
  & = & \{ X \in {}_{ Q,A }\!\Mod \mid \mbox{$\BH_1^{ [q] }( X ) = 0$ for each $q \in Q_0$} \}, \\
  \end{array}
$
\end{itemize}
see \cite[thm.\ 7.1]{HJ-JLMS}.  By the same theorem we have
\begin{itemize}
\setlength\itemsep{4pt}

  \item
\vspace{-4ex}
$
  \begin{array}{rcl}
  \phantom{x} && \\[2mm]
  \cE & = & \{ X \in {}_{ Q,A }\!\Mod \mid \mbox{$\BH^i_{ [q] }( X ) = 0$ for each $i \in \BZ$ and each $q \in Q_0$} \} \\[2mm]
  & = & \{ X \in {}_{ Q,A }\!\Mod \mid \mbox{$\BH_i^{ [q] }( X ) = 0$ for each $i \in \BZ$ and each $q \in Q_0$} \}. \\
  \end{array}
$
\end{itemize}
Combining the last bullet with the long exact $\Ext$ sequence for $\BH^i_{ [q] }$ and the long exact $\Tor$ sequence for $\BH_i^{ [q] }$ makes it easy to see that
\begin{itemize}
\setlength\itemsep{4pt}
\setcounter{enumi}{2}

  \item  $\cE$ is a wide subcategory in the sense of \ref{bfhpg:wide} in Appendix \ref{app:II}.

\end{itemize}
\end{bfhpg}

\begin{bfhpg}[The class $\weq$ of weak equivalences]
\label{bfhpg:weq}
The class of {\em weak equivalences} in ${}_{ Q,A }\Mod$ is
\begin{itemize}
\setlength\itemsep{4pt}
\setcounter{enumi}{1}

  \item
\vspace{-9ex}
$
  \begin{array}{rcl}
  \phantom{x} && \\[11mm]  
  \weq & = &
    \Bigg\{
    \mbox{morphisms $\varphi$ in ${}_{ Q,A }\!\Mod$}\: 
    \Bigg|
    \begin{array}{l}
      \mbox{$\BH^1_{ [q] }( \varphi )$ and $\BH^2_{ [q] }( \varphi )$ are}\\[1mm]
      \mbox{isomorphisms for each $q \in Q_0$} 
    \end{array}
    \Bigg\} \\[7mm]
  & = &
    \Bigg\{
    \mbox{morphisms $\varphi$ in ${}_{ Q,A }\!\Mod$}\:
    \Bigg|
    \begin{array}{l}
      \mbox{$\BH_1^{ [q] }( \varphi )$ and $\BH_2^{ [q] }( \varphi )$ are}\\[1mm]
      \mbox{isomorphisms for each $q \in Q_0$} 
    \end{array}
    \Bigg\}, \\
  \end{array}
$
\end{itemize}
see \cite[thm.\ 7.2]{HJ-JLMS}.  By the same theorem we have
\begin{itemize}
\setlength\itemsep{4pt}

  \item
\vspace{-9ex}
$
  \begin{array}{rcl}
  \phantom{x} && \\[11mm]  
  \weq & = &
    \Bigg\{
    \mbox{morphisms $\varphi$ in ${}_{ Q,A }\!\Mod$}\: 
    \Bigg|
    \begin{array}{l}
      \mbox{$\BH^i_{ [q] }( \varphi )$ is an isomorphism for}\\[1mm]
      \mbox{each $i \in \BZ$ and each $q \in Q_0$} 
    \end{array}
    \Bigg\} \\[7mm]
  & = &
    \Bigg\{
    \mbox{morphisms $\varphi$ in ${}_{ Q,A }\!\Mod$}\:
    \Bigg|
    \begin{array}{l}
      \mbox{$\BH_i^{ [q] }( \varphi )$ is an isomorphism for}\\[1mm]
      \mbox{each $i \in \BZ$ and each $q \in Q_0$} 
    \end{array}
    \Bigg\}. \\
  \end{array}
$
\end{itemize}
Perhaps surprisingly, weak equivalences are not in general characterised by $\BH^1_{ [q] }( \varphi )$ alone being an isomorphism for each $q \in Q_0$ or by $\BH_1^{ [q] }( \varphi )$ alone being an isomorphism for each $q \in Q_0$; see \cite[exa.\ 8.21]{HJ-JLMS}.
\end{bfhpg}

\begin{bfhpg}[The functors $\BH_{ [q] }^i$ and $\BH^{ [q] }_i$, the class $\cE$, and the class $\weq$ for complexes]
\label{bfhpg:E_etc_for_complexes}
Let $Q = \Qcpx$ whence ${}_{ Q,A }\!\Mod = \Ch( A )$; see \ref{bfhpg:examples_of_Q}.  An object $X \in {}_{ Q,A }\!\Mod$ is a complex
$
  X = \cdots
      \xrightarrow{} X_2
      \xrightarrow{} X_1
      \xrightarrow{} X_0
      \xrightarrow{} X_{-1}
      \xrightarrow{} X_{-2}
      \xrightarrow{} \cdots
$
over $\Mod( A )$ and
\begin{itemize}
\setlength\itemsep{4pt}

  \item  $\BH^i_{ [q] }( X ) = H_{ q-i }( X ),$

  \item  $\BH_i^{ [q] }( X ) = H_{ q+i }( X )$

\end{itemize}
for $i \geqslant 1$ where $H_j$ is classic homology at degree $j$.  It follows that
\begin{itemize}
\setlength\itemsep{4pt}
\setcounter{enumi}{2}

  \item  $\cE = \{ \mbox{ exact complexes } \}$,
  
  \item  $\weq = \{ \mbox{ quasiisomorphisms } \}$.

\end{itemize}
Note that in this particular case, weak equivalences are in fact characterised by $\BH^1_{ [q] }( \varphi )$ being an isomorphism for each $q \in Q_0$ and by $\BH_1^{ [q] }( \varphi )$ being an isomorphism for each $q \in Q_0$.
\end{bfhpg}

\begin{bfhpg}[The functors $\BH_{ [q] }^i$ and $\BH^{ [q] }_i$, the class $\cE$, and the class $\weq$ for $N$-complexes]
\label{bfhpg:E_etc_for_N-complexes}
Let $Q = \QNcpx$ whence ${}_{ Q,A }\!\Mod = \ChN( A )$; see \ref{bfhpg:examples_of_Q}.  An object $X \in {}_{ Q,A }\!\Mod$ is an $N$-complex
$
  X = \cdots
      \xrightarrow{} X_2
      \xrightarrow{} X_1
      \xrightarrow{} X_0
      \xrightarrow{} X_{-1}
      \xrightarrow{} X_{-2}
      \xrightarrow{} \cdots
$
over $\Mod( A )$.

Given $q \in Q_0$ and an integer $0 < j < N$, there is a generalised homology group given by
\[
  {}_{ j }H_q( X ) =
  \frac{ \Ker( X_q \xrightarrow{} \cdots \xrightarrow{} X_{ q-j } ) }{ \Image( X_{ q+N-j } \xrightarrow{} \cdots \xrightarrow{} X_q ) },
\]
see \cite[def.\ 1.1]{Kapranov}.  The definition makes sense because $\Image \subseteq \Ker$ since the composition $X_{ q+N-j } \xrightarrow{} \cdots \xrightarrow{} X_{ q-j }$ is $0$ in the $N$-complex  $X$.  In these terms, for $i \geqslant 1$,
\begin{itemize}
\setlength\itemsep{6pt}

  \item  $\BH^i_{ [q] }( X ) =
\left\{
  \begin{array}{cl}
    {}_{ N-1 \vphantom{q-1-\frac{i-1}{2}N} }H_{ q-1-\frac{i-1}{2}N }( X ) & \mbox{ for $i$ odd, } \\[3mm]
    {}_{ 1 \vphantom{q-\frac{i}{2}N} }H_{ q-\frac{i}{2}N }( X ) & \mbox{ for $i$ even, } \\
  \end{array}
\right.  
$

  \item  $\BH_i^{ [q] }( X ) = 
\left\{
  \begin{array}{cl}
    {}_{ N-1 \vphantom{q+1+\frac{i-1}{2}N} }H_{ q+1+\frac{i-1}{2}N }( X ) & \mbox{ for $i$ odd, } \\[3mm]
    {}_{ 1 \vphantom{q+\frac{i}{2}} }H_{ q+\frac{i}{2}N }( X ) & \mbox{ for $i$ even. } \\
  \end{array}
\right.  
$

\end{itemize}

It follows that
\begin{itemize}
\setlength\itemsep{4pt}
\setcounter{enumi}{2}

  \item  $\cE = \{ \mbox{ $N$-exact $N$-complexes } \}$; see \cite[def.\ 1.1 and prop.\ 1.5]{Kapranov}.
  
  \item  $\weq = \{ \mbox{ $N$-quasiisomorphisms } \}$; see \cite[def.\ 3.6]{Iyama-Kato-Miyachi-N} and \cite[def.\ 1.1 and prop.\ 1.5]{Kapranov}.

\end{itemize}
Note that an $N$-complex $X$ is called {\rm $N$-exact} if ${}_{ j }H_q( X ) = 0$ for each $q \in Q_0$ and each integer $0 < j < N$.  A morphism $X \xrightarrow{ \varphi } Y$ of $N$-complexes is called an {\em $N$-quasiisomorphism} if ${}_{ j }H_q( \varphi )$ is an isomorphism for each $q \in Q_0$ and each integer $0 < j < N$.   
\end{bfhpg}

\section{The Frobenius categories ${}^{ \perp }\cE$ and $\cE^{ \perp }$}
\label{sec:E_perp}

\begin{bfhpg}[The Frobenius categories ${}^{ \perp }\cE$ and $\cE^{ \perp }$]
\label{bfhpg:Frobenius_E}
The full subcategory $\cE$ of exact objects was introduced in \ref{bfhpg:E} and the notation $\perp$ in \ref{bfhpg:QAMod}.  They permit the following definitions.
\begin{itemize}
\setlength\itemsep{4pt}

  \item  The full subcategory of {\em semiprojective} objects is ${}^{ \perp }\cE$.

  \item  The full subcategory of {\em semiinjective} objects is $\cE^{ \perp }$.

\end{itemize}
These are all full subcategories of ${}_{ Q,A }\!\Mod$.  Since ${}^{ \perp }\cE$ is defined as an $\Ext^1_{ Q,A }$ perpendicular subcategory, it is closed under extensions in ${}_{ Q,A }\!\Mod$ and satisfies ${}_{ Q,A }\!\Prj \subseteq {}^{ \perp }\cE$.  Extension closure implies that ${}^{ \perp }\cE$ is an exact category in a canonical way, see \ref{bfhpg:exact}.  The conflations are the short exact sequences in the abelian category ${}_{ Q,A }\!\Mod$ which consist of objects from ${}^{ \perp }\cE$.  

By \cite[thm.\ 6.5]{HJ-JLMS} these statements can be strengthened to the first of the following, and the second is analogous.  The notion of Frobenius category is explained in \ref{bfhpg:Frobenius}.
\begin{itemize}
\setlength\itemsep{4pt}
\setcounter{enumi}{2}

  \item  ${}^{ \perp }\cE$ is a Frobenius category with ${}_{ Q,A }\!\Prj$ as its class of projective-injective objects. 
	
  \item  $\cE^{ \perp }$ is a Frobenius category with ${}_{ Q,A }\!\Inj$ as its class of projective-injective objects.

\end{itemize}
This is proved using the theory of model categories to which we shall return in Section \ref{sec:model}.

Giving a precise description of the objects in ${}^{ \perp }\cE$ and $\cE^{ \perp }$ is in general difficult, but by \cite[thm.\ E]{HJ-arXiv} we do have the following inclusions which are equalities if the left global dimension of $A$ is finite.
\begin{itemize}
\setlength\itemsep{4pt}
\setcounter{enumi}{4}

  \item  ${}^{ \perp }\cE \subseteq \{ X \in {}_{ Q,A }\Mod \mid \mbox{$X( q )$ is projective for each $q \in Q_0$} \}$ 

  \item  $\cE^{ \perp } \subseteq \{ X \in {}_{ Q,A }\Mod \mid \mbox{$X( q )$ is injective for each $q \in Q_0$} \}$ 

\end{itemize}
\end{bfhpg}

\begin{bfhpg}[The Frobenius categories ${}^{ \perp }\cE$ and $\cE^{ \perp }$ for complexes]
\label{bfhpg:Eperp_for_complexes}
Let $Q = \Qcpx$ whence ${}_{ Q,A }\!\Mod = \Ch( A )$.  Then
\begin{itemize}
\setlength\itemsep{4pt}

  \item  We have
\[  
  {}^{ \perp }\cE =
  \{ \mbox{ semiprojective complexes } \}
\]  
by \cite[prop.\ 2.3.5]{Garcia-Rozas-book} (which uses the term ``DG-projective'').  Semiprojective complexes were introduced in \cite[sec.\ 2]{Boekstedt-Neeman} under the name ``special complexes of projectives''.  They consist of projective modules by \ref{bfhpg:Frobenius_E} and include all right bounded complexes of projective modules, in particular the projective resolution of each $A$-left module.
  
  \item  We have
\[  
  \cE^{ \perp } = 
  \{ \mbox{ semiinjective complexes } \}
\]  
by \cite[prop.\ 2.3.4]{Garcia-Rozas-book} (which uses the term ``DG-injective'').  Semiinjective complexes were also introduced in \cite[sec.\ 2]{Boekstedt-Neeman} under the name ``special complexes of injectives''.  They consist of injective modules by \ref{bfhpg:Frobenius_E} and include all left bounded complexes of injective modules, in particular the injective resolution of each $A$-left module.

\end{itemize}
\end{bfhpg}

\begin{bfhpg}[The Frobenius categories ${}^{ \perp }\cE$ and $\cE^{ \perp }$ for $N$-complexes]
\label{bfhpg:Eperp_for_N-complexes}
Let $Q = \QNcpx$ whence ${}_{ Q,A }\!\Mod = \ChN( A )$.  Then
\begin{itemize}
\setlength\itemsep{4pt}

  \item  The objects in ${}^{ \perp }\cE$ will be called semiprojective $N$-complexes.  These $N$-complexes were considered in \cite[ex.\ 3.6]{Bahiraei} under the name ``dg-projective $N$-complexes''.  They consist of projective modules by \ref{bfhpg:Frobenius_E} and include all right bounded $N$-complexes of projective modules by the lines immediately after \cite[def.\ 3.16]{Iyama-Kato-Miyachi-N}.  They are included in the $K$-projective $N$-complexes of \cite[def.\ 3.16]{Iyama-Kato-Miyachi-N}.

  \item  The objects in $\cE^{ \perp }$ will be called semiinjective $N$-complexes.  These $N$-complexes consist of injective modules by \ref{bfhpg:Frobenius_E} and include all left bounded $N$-complexes of injective modules by the lines immediately after \cite[def.\ 3.16]{Iyama-Kato-Miyachi-N}.  They are included in the $K$-injective $N$-complexes of \cite[def.\ 3.16]{Iyama-Kato-Miyachi-N}.  

\end{itemize}
\end{bfhpg}

\section{The $Q$-shaped derived category $\cD_Q( A )$}
\label{sec:DQ}

\begin{bfhpg}[The $Q$-shaped derived category $\cD_Q( A )$]
\label{bfhpg:DQA}
The {\em $Q$-shaped derived category of $A$} is the homotopy category
\[
  \cD_Q( A ) 
  = \Ho( {}_{ Q,A }\!\Mod )
  = \weq^{ -1 }\!{}_{ Q,A }\!\Mod
\]
obtained from ${}_{ Q,A }\!\Mod$ by formally inverting each weak equivalence.  There are equivalences of categories
\[
  \frac{ {}^{ \perp }\cE }{ {}_{ Q,A }\Prj }
  \cong \cD_Q( A ) \cong
  \frac{ \cE^{ \perp } }{ {}_{ Q,A }\Inj }.
\]
Here $\frac{ {}^{ \perp }\cE }{ {}_{ Q,A }\Prj }$ and $\frac{ \cE^{ \perp } }{ {}_{ Q,A }\Inj }$ are the stable categories of the Frobenius categories ${}^{ \perp }\cE$ and $\cE^{ \perp }$; see \ref{bfhpg:stable}.  Hence they are triangulated categories, and in this part of the paper, we view them as the de facto definition of $\cD_Q( A )$.

The equivalences are established using the theory of model categories; see Section \ref{sec:model}.
\end{bfhpg}

\begin{bfhpg}[The $Q$-shaped derived category $\cD_Q( A )$ for complexes]
\label{bfhpg:DQA_cpx}
Let $Q = \Qcpx$ whence ${}_{ Q,A }\!\Mod = \Ch( A )$.  There is an equivalence of triangulated categories
\[
  \cD_Q( A ) \cong \cD( A ),
\]  
where $\cD( A )$ is the classic derived category.  To see this, note that the objects in $\cD_Q( A ) \cong \frac{ {}^{ \perp }\cE }{ {}_{ Q,A }\!\Prj }$ are the semiprojective complexes by \ref{bfhpg:Eperp_for_complexes}.  The morphisms are chain maps modulo chain maps which factorise through a projective object.  Such factorisations exist precisely for null homotopic chain maps, as one can see by amending the arguments in \cite[p.\ 28]{Happel-book}.  Hence $\frac{ {}^{ \perp }\cE }{ {}_{ Q,A }\Prj }$ is the category of semiprojective complexes and chain maps modulo chain homotopy, which by \cite[p.\ 216]{Boekstedt-Neeman} is triangulated equivalent to $\cD( A )$.  
\end{bfhpg}

\begin{bfhpg}[The $Q$-shaped derived category $\cD_Q( A )$ for $N$-complexes]
\label{bfhpg:DQA_Ncpx}
Let $Q = \QNcpx$ whence ${}_{ Q,A }\!\Mod = \ChN( A )$.  There is an equivalence of triangulated categories
\[
  \cD_Q( A ) \cong \cD_N( A ),
\]  
where $\cD_N( A )$ is the derived category of $N$-complexes of \cite[def.\ 3.6]{Iyama-Kato-Miyachi-N}.  To see this, note that the objects in $\cD_Q( A ) \cong \frac{ {}^{ \perp }\cE }{ {}_{ Q,A }\Prj }$ are the semiprojective $N$-complexes from \ref{bfhpg:Eperp_for_N-complexes}.  The morphisms are morphisms of $N$-complexes modulo morphisms which factorise through a projective object.  Such factorisations exist precisely for morphisms which are $N$-null homotopic, see \cite[p.\ 8]{Kapranov}, as one can see by amending the proof of \cite[thm.\ 2.3]{Iyama-Kato-Miyachi-N}.  Hence $\frac{ {}^{ \perp }\cE }{ {}_{ Q,A }\Prj }$ is the category of semiprojective $N$-complexes and morphisms of $N$-complexes modulo $N$-chain homotopy, which by \cite[thm.\ 3.17(i) and its proof]{Iyama-Kato-Miyachi-N} is triangulated equivalent to $\cD_N( A )$.  Note that the proof rather than the formulation of \cite[thm.\ 3.17(i)]{Iyama-Kato-Miyachi-N} shows that the semiprojective $N$-complexes suffice in this statement.
\end{bfhpg}

\section{The suspension functor of $\cD_Q( A )$}
\label{sec:suspension}

The stable category of a Frobenius category is triangulated.  In particular, it has a suspension functor which can be computed as described in the last bullet point of \ref{bfhpg:stable}.  We will do the computation for $\cD_Q( A )$ in our two standing examples, complexes and $N$-complexes, and also in the case of $m$-periodic complexes.

\begin{bfhpg}[The suspension functor for complexes]
\label{bfhpg:suspension_for_complexes}
Let $Q = \Qcpx$ whence ${}_{ Q,A }\!\Mod = \Ch( A )$.  Consider the Frobenius category ${}^{ \perp }\cE$ and its stable category $\frac{ {}^{ \perp }\cE }{ {}_{ Q,A }\Prj }$.  We will compute the action of the suspension functor $\Sigma$ on objects using \ref{bfhpg:stable}.  The computation will show that $\Sigma$ acts as one would expect from classic homological algebra: It shifts a complex one step against the direction of the differential and flips the sign of the differential.

Let $P \in {}^{ \perp }\cE$ be given.  The following diagram shows a conflation $P \xrightarrow{} R \xrightarrow{} P'$ in ${}^{ \perp }\cE$ with $R \in {}_{ Q,A }\!\Prj$ whence $\Sigma P = P'$.
\[
\vcenter{
\xymatrix @=3.5pc
{
  P=\cdots \ar[r] \ar@<-2.9ex>[d] & P_{i+1} \ar^{\partial_{i+1}}[r] \ar_{ \begin{psmallmatrix} \id \\ \partial_{i+1} \end{psmallmatrix} }[d] & P_i \ar^{\partial_i}[r] \ar^{ \begin{psmallmatrix} \id \\ \partial_i \end{psmallmatrix} }[d] & P_{i-1} \ar[r] \ar^{ \begin{psmallmatrix} \id \\ \partial_{i-1} \end{psmallmatrix} }[d] & \cdots \\
  R=\cdots \ar[r] \ar@<-2.9ex>[d] & P_{i+1} \oplus P_i \ar[r]^{ \begin{psmallmatrix} 0 & \id \\ 0 & 0 \end{psmallmatrix} } \ar_{ \begin{psmallmatrix} -\partial_{i+1} & \id \end{psmallmatrix} }[d] & P_i \oplus P_{i-1} \ar[r]^{ \begin{psmallmatrix} 0 & \id \\ 0 & 0 \end{psmallmatrix} } \ar^{ \begin{psmallmatrix} -\partial_i & \id \end{psmallmatrix} }[d] & P_{i-1} \oplus P_{i-2} \ar[r] \ar^{ \begin{psmallmatrix} -\partial_{ i-1 } & \id \end{psmallmatrix} }[d] & \cdots \\
  P'=\cdots \ar[r] 
   & P_i \ar_{-\partial_i}[r] 
   & P_{i-1} \ar_{-\partial_{i-1}}[r] 
   & P_{i-2} \ar[r] 
   & \cdots \\
}
        }
\]
To see that $R \in {}_{ Q,A }\!\Prj$, note that $R$ is a direct sum of shifts of complexes of the form
\[
  \cdots \xrightarrow{} 0 \xrightarrow{} P \xrightarrow{ \id } P \xrightarrow{} 0 \xrightarrow{} \cdots
\]
with $P$ a projective $A$-module, and each such is easily checked to be in ${}_{ Q,A }\Prj$.  

To see that the diagram shows a conflation in ${}^{ \perp }\cE$, we must first see that it shows a short exact sequence in ${}_{ Q,A }\!\Mod$, and this is true because it is split exact in each degree.  Secondly, we must see that $P, R, P' \in {}^{ \perp }\cE$.  This holds for $P$ by definition, for $P'$ because $P'$ is isomorphic to a shift of $P$, and for $R$ because $R \in {}_{ Q,A }\!\Prj$.
\end{bfhpg}

\begin{bfhpg}[The suspension functor for $N$-complexes]
\label{bfhpg:suspension_for_N-complexes}
Let $Q = \QNcpx$ whence ${}_{ Q,A }\!\Mod = \ChN( A )$.  Consider the Frobenius category ${}^{ \perp }\cE$ and its stable category $\frac{ {}^{ \perp }\cE }{ {}_{ Q,A }\Prj }$.  The action of the suspension functor $\Sigma$ on objects is computed in \cite[p.\ 693]{Iyama-Kato-Miyachi-N} using the method of \ref{bfhpg:stable}.  We will not show the full computation but give an example to illustrate that $\Sigma$ is not given by shifting.

Let $P$ be a projective $A$-module and consider $P$ as an $N$-complex concentrated in degree $0$.  The following diagram shows a conflation $P \xrightarrow{} R \xrightarrow{} P'$ in ${}^{ \perp }\cE$ with $R \in {}_{ Q,A }\!\Prj$ whence $\Sigma P = P'$.
\begingroup
\[
\arraycolsep=1.85pt
  \begin{array}{ccccccccccccccccccc}
    &&&&&&&&&&&&&& \hspace{-6ex} \lefteqn{\mbox{\tiny degree $0$}} &&&& \\
    &&&&&&&&&&&&&& \downarrow &&&& \\
    P
    & \hspace{0.31ex} =
    & \cdots
    & \xrightarrow{\phantom{1_P}}
    & \mathmakebox[\widthof{$\cdots$}][c]{0}
    & \xrightarrow{\phantom{1_P}}    
    & \mathmakebox[\widthof{$\cdots$}][c]{0}
    & \xrightarrow{\phantom{1_P}}
    & \mathmakebox[\widthof{$\cdots$}][c]{0}
    & \xrightarrow{\phantom{1_P}}
    & \cdots
    & \xrightarrow{\phantom{1_P}}
    & \mathmakebox[\widthof{$\cdots$}][c]{0}
    & \xrightarrow{\phantom{1_P}}
    & \mathmakebox[\widthof{$\cdots$}][c]{P}
    & \xrightarrow{\phantom{1_P}}
    & \mathmakebox[\widthof{$\cdots$}][c]{0}
    & \xrightarrow{\phantom{1_P}}
    & \cdots \\
    \downarrow &&&& \downarrow && \downarrow && \downarrow &&&& \downarrow && \,\,\downarrow \lefteqn{ \scriptstyle 1_P } && \downarrow \\
    R
    & \hspace{0.31ex} =
    & \cdots
    & \xrightarrow{\phantom{1_P}}
    & \mathmakebox[\widthof{$\cdots$}][c]{0}
    & \xrightarrow{\phantom{1_P}}    
    & \mathmakebox[\widthof{$\cdots$}][c]{P}
    & \xrightarrow{1_P}
    & \mathmakebox[\widthof{$\cdots$}][c]{P}
    & \xrightarrow{1_P}
    & \cdots
    & \xrightarrow{1_P}
    & \mathmakebox[\widthof{$\cdots$}][c]{P}
    & \xrightarrow{1_P}
    & \mathmakebox[\widthof{$\cdots$}][c]{P}
    & \xrightarrow{\phantom{1_P}}
    & \mathmakebox[\widthof{$\cdots$}][c]{0}
    & \xrightarrow{\phantom{1_P}}
    & \cdots \\    
    \downarrow &&&& \downarrow && \,\,\downarrow \lefteqn{ \scriptstyle 1_P } && \,\,\downarrow \lefteqn{ \scriptstyle 1_P } &&&& \,\,\downarrow \lefteqn{ \scriptstyle 1_P } && \downarrow && \downarrow \\    
    P'
    & \hspace{0.31ex} =
    & \cdots
    & \xrightarrow{\phantom{1_P}}
    & \mathmakebox[\widthof{$\cdots$}][c]{0}
    & \xrightarrow{\phantom{1_P}}    
    & \mathmakebox[\widthof{$\cdots$}][c]{P}
    & \xrightarrow[1_P]
    & \mathmakebox[\widthof{$\cdots$}][c]{P}
    & \xrightarrow[1_P]
    & \cdots
    & \xrightarrow[1_P]
    & \mathmakebox[\widthof{$\cdots$}][c]{P}
    & \xrightarrow{\phantom{1_P}}
    & \mathmakebox[\widthof{$\cdots$}][c]{0}
    & \xrightarrow{\phantom{1_P}}
    & \mathmakebox[\widthof{$\cdots$}][c]{0}
    & \xrightarrow{\phantom{1_P}}
    & \cdots \\    
    &&&&&& \uparrow &&&& \\
    &&&&&& \hspace{-9ex} \lefteqn{\mbox{\tiny degree $N-1$}} &&&& \\
  \end{array}
\]
\endgroup
It is easy to check $R \in {}_{ Q,A }\!\Prj$.  To see that the diagram shows a conflation in ${}^{ \perp }\cE$, note that it clearly shows a short exact sequence and that $P, R, P' \in {}^{ \perp }\cE$ because $P, R, P'$ are bounded complexes of projective modules; see \ref{bfhpg:Eperp_for_N-complexes}.

As shown by this example, there are categories $Q$ which give considerably more complicated formulae for $\Sigma$ than the ``shift plus sign flip'' that applies to complexes in classic homological algebra.  
\end{bfhpg}

\begin{bfhpg}[The suspension functor for $m$-periodic complexes]
\label{bfhpg:suspension_for_m-periodic_complexes}
Let $m \geqslant 1$ be an integer and let $Q$ be defined by Figure \ref{fig:cyclic_quiver} modulo the relations that consecutive arrows compose to $0$.  Then ${}_{ Q,A }\!\Mod$ is the category of $m$-periodic chain complexes and chain maps.  Consider the Frobenius category ${}^{ \perp }\cE$ and its stable category $\frac{ {}^{ \perp }\cE }{ {}_{ Q,A }\Prj }$.

The action of the suspension functor $\Sigma$ on objects can be computed by the same argument as in \ref{bfhpg:suspension_for_complexes}, and the conclusion is the same: $\Sigma$ shifts a complex one step against the direction of the differential and flips the sign of the differential.  A slightly more elaborate argument shows that $\Sigma$ also shifts morphisms one step against the direction of the differential.

It follows by $m$-periodicity that if $m$ is even then $\Sigma^m \cong \id$.  If $m$ is odd, then $\Sigma^{2m} \cong \id$, where the additional factor $2$ is necessary because $\Sigma^m$ flips the sign of the differential.
\end{bfhpg}

\part{The model category approach to $\cD_Q( A )$}
\label{part:model}

This part constructs the projective and injective model category structures on ${}_{ Q,A }\!\Mod$ and obtains $\cD_Q( A )$ as the corresponding homotopy category where each weak equivalence has been formally inverted.

\section{Cotorsion pairs in ${}_{ Q,A }\!\Mod$}
\label{sec:cotorsion}

This section introduces four cotorsion pairs which are hereditary and functorially complete.  See \ref{bfhpg:cotorsion_pair} for the relevant definitions.

\begin{bfhpg}[The cotorsion pair $( {}^{ \perp }\cE,\cE )$]
\label{bfhpg:cotorsion_pair1}
The class $\cE$ of exact objects and its $\Ext^1$-perpendicular ${}^{ \perp }\cE$ were introduced in \ref{bfhpg:E} and \ref{bfhpg:Frobenius_E}.  They are both full subcategories of ${}_{ Q,A }\!\Mod$.
\begin{itemize}
\setlength\itemsep{4pt}

  \item  $( {}^{ \perp }\cE,\cE )$ is a cotorsion pair in ${}_{ Q,A }\!\Mod$.  It satisfies ${}^{ \perp }\cE \cap \cE = {}_{ Q,A }\!\Prj$ and it is hereditary and  functorially complete.

\end{itemize}
These claims are proved in \cite[thm.\ 4.4(a)]{HJ-JLMS} except functorial completeness.  To prove this as well, observe that \cite[proof of thm. 5.5]{HJ-JLMS} produces a set $\cS$ of objects of ${}_{ Q,A }\!\Mod$ such that $\cS^{ \perp } = \cE$.  This property is preserved by adding to $\cS$ a projective generator of ${}_{ Q,A }\!\Mod$, which exists by \ref{bfhpg:QAMod}.  Hence $( {}^{ \perp }\cE,\cE )$ is functorially complete by \cite[thm.\ 2.1]{Gillespie-BLMS}.  See also \cite[thm.\ 5.16]{Stovicek-exact}.
\end{bfhpg}

\begin{bfhpg}[The cotorsion pair $( {}_{ Q,A }\!\Prj,{}_{ Q,A }\!\Mod )$]
\label{bfhpg:cotorsion_pair2}
Consider the category ${}_{ Q,A }\!\Mod$ and its full subcategory ${}_{ Q,A }\!\Prj$ of projective objects.
\begin{itemize}
\setlength\itemsep{4pt}

  \item  $( {}_{ Q,A }\!\Prj,{}_{ Q,A }\!\Mod )$ is a cotorsion pair in ${}_{ Q,A }\!\Mod$.  It is hereditary and functorially complete.

\end{itemize}
The hereditary property is immediate.  To prove the remaining claims, observe that ${}_{ Q,A }\!\Mod$ has a projective generator $P$ by \ref{bfhpg:QAMod}.  Hence $( {}_{ Q,A }\!\Prj,{}_{ Q,A }\!\Mod ) = \big( {}^{ \perp }( \{ P \}^{ \perp } ),\{ P \}^{ \perp } \big)$ is a functorially complete cotorsion pair by \cite[thm.\ 2.1]{Gillespie-BLMS}.
\end{bfhpg}

\begin{bfhpg}[The cotorsion pair $( \cE,\cE^{ \perp } )$]
\label{bfhpg:cotorsion_pair3}
The class $\cE$ of exact objects and its $\Ext^1$-perpendicular $\cE^{ \perp }$ were introduced in \ref{bfhpg:E} and \ref{bfhpg:Frobenius_E}.  They are both full subcategories of ${}_{ Q,A }\!\Mod$.
\begin{itemize}
\setlength\itemsep{4pt}

  \item  $( \cE,\cE^{ \perp } )$ is a cotorsion pair in ${}_{ Q,A }\!\Mod$.  It satisfies $\cE \cap \cE^{ \perp } = {}_{ Q,A }\!\Inj$ and it is hereditary and functorially complete.

\end{itemize}
These claims are proved in \cite[thm.\ 4.4(b)]{HJ-JLMS} except functorial completeness.  To prove this as well, observe that \cite[proof of thm.\ 5.9]{HJ-JLMS} verifies that $\cE$ satisfies conditions (1) and (2) of \cite[thm.\ A.3]{HJ-JLMS}.  Hence \cite[proof of thm.\ A.3]{HJ-JLMS} implies that \cite[cor.\ 5.17]{Stovicek-exact} applies to $\cE$ whence $( \cE,\cE^{ \perp } )$ is functorially complete.
\end{bfhpg}

\begin{bfhpg}[The cotorsion pair $( {}_{ Q,A }\!\Mod,{}_{ Q,A }\!\Inj )$]
\label{bfhpg:cotorsion_pair4}
Consider the category ${}_{ Q,A }\!\Mod$ and its full subcategory ${}_{ Q,A }\!\Inj$ of injective objects.
\begin{itemize}
\setlength\itemsep{4pt}

  \item  $( {}_{ Q,A }\!\Mod,{}_{ Q,A }\!\Inj )$ is a cotorsion pair in ${}_{ Q,A }\!\Mod$.  It is hereditary and functorially complete.

\end{itemize}
The hereditary property is immediate.  The remaining claims are proved in \cite[cor.\ 5.9]{Stovicek-exact}.  To see that this result applies to ${}_{ Q,A }\!\Mod$, observe that this category is, in the terminology of \cite[cor.\ 5.9]{Stovicek-exact}, exact of Grothendieck type.  This holds by \cite[text after def.\ 3.11 and prop.\ 3.13]{Stovicek-exact}.
\end{bfhpg}

\section{The projective and injective model category structures on ${}_{ Q,A }\!\Mod$}
\label{sec:model}

This section introduces two hereditary Hovey triples.  It proceeds to study the ensuing so-called projective and injective model category structures on ${}_{ Q,A }\!\Mod$.  See \ref{bfhpg:model_categories}, \ref{bfhpg:abelian_model_categories}, and \ref{bfhpg:Hovey_triple} for the relevant definitions.

\begin{bfhpg}[The Hovey triple $( {}^{ \perp }\cE,\cE,{}_{ Q,A }\!\Mod )$]
\label{bfhpg:Hovey_triple_P}
In the abelian category ${}_{ Q,A }\!\Mod$,
\begin{itemize}
\setlength\itemsep{4pt}

  \item  $( \cC_{ \psub },\cW_{ \psub },\cF_{ \psub } ) = ( {}^{ \perp }\cE,\cE,{}_{ Q,A }\!\Mod )$ is a hereditary Hovey triple.

\end{itemize}
To see this, we check the conditions in \ref{bfhpg:Hovey_triple}.
\begin{itemize}
\setlength\itemsep{4pt}
\setcounter{enumi}{1}

  \item  $\cW_{ \psub } = \cE$ is wide by \ref{bfhpg:E}.
  
  \item  We have $( \cC_{ \psub },\cW_{ \psub } \cap \cF_{ \psub } ) = ( {}^{ \perp }\cE,\cE \cap {}_{ Q,A }\!\Mod ) = ( {}^{ \perp }\cE,\cE )$.  This is a hereditary functorially complete cotorsion pair by \ref{bfhpg:cotorsion_pair1}. 
  
  \item  We have $( \cC_{ \psub } \cap \cW_{ \psub },\cF_{ \psub } ) = ( {}^{ \perp }\cE \cap \cE,{}_{ Q,A }\!\Mod ) = ( {}_{ Q,A }\!\Prj,{}_{ Q,A }\!\Mod )$ by \ref{bfhpg:cotorsion_pair1}.  This is a hereditary functorially complete cotorsion pair by \ref{bfhpg:cotorsion_pair2}. 

\end{itemize}
\end{bfhpg}

\begin{bfhpg}[The Hovey triple $( {}_{ Q,A }\!\Mod,\cE,\cE^{ \perp } )$]
\label{bfhpg:Hovey_triple_I}
In the abelian category ${}_{ Q,A }\!\Mod$,
\begin{itemize}
\setlength\itemsep{4pt}

  \item  $( \cC_{ \isub },\cW_{ \isub },\cF_{ \isub } ) = ( {}_{ Q,A }\!\Mod,\cE,\cE^{ \perp } )$ is a hereditary Hovey triple.

\end{itemize}
To see this, we check the conditions in \ref{bfhpg:Hovey_triple}.
\begin{itemize}
\setlength\itemsep{4pt}
\setcounter{enumi}{1}

  \item  $\cW_{ \isub } = \cE$ is wide by \ref{bfhpg:E}.
  
  \item  We have $( \cC_{ \isub },\cW_{ \isub } \cap \cF_{ \isub } ) = ( {}_{ Q,A }\!\Mod,\cE \cap \cE^{ \perp } ) = ( {}_{ Q,A }\!\Mod,{}_{ Q,A }\!\Inj )$ by \ref{bfhpg:cotorsion_pair3}.  This is a hereditary functorially complete cotorsion pair by  \ref{bfhpg:cotorsion_pair4}. 
    
  \item  We have $( \cC_{ \isub } \cap \cW_{ \isub },\cF_{ \isub } ) = ( {}_{ Q,A }\!\Mod \cap\, \cE,\cE^{ \perp } ) = ( \cE,\cE^{ \perp } )$.  This is a hereditary functorially complete cotorsion pair by \ref{bfhpg:cotorsion_pair3}. 

\end{itemize}
\end{bfhpg}

\begin{bfhpg}[The projective model category structure on ${}_{ Q,A }\!\Mod$]
\label{bfhpg:projective_model_category_structure}
By Theorem \ref{thm:Hovey}, the hereditary Hovey triple $( \cC_{ \psub },\cW_{ \psub },\cF_{ \psub } ) = ( {}^{ \perp }\cE,\cE,{}_{ Q,A }\!\Mod )$ from \ref{bfhpg:Hovey_triple_P} gives a model category structure $( \weq_{ \psub },\cof_{ \psub },\fib_{ \psub } )$ on ${}_{ Q,A }\!\Mod$, called the {\em projective model category structure}, which can be described as follows.
\begin{itemize}
\setlength\itemsep{4pt}

  \item  $\weq_{ \psub }$ consists of the compositions $\pi\iota$ where $\pi$ is an epimorphism with kernel in $\cW_{ \psub } = \cE$ and $\iota$ is a monomorphism with cokernel in $\cW_{ \psub } = \cE$.  Note that
\[  
  \weq_{ \psub } = \weq
\]
by \cite[thm.\ 7.2]{HJ-JLMS}, where $\weq$ is the class from \ref{bfhpg:weq}.

  \item  $\cof_{ \psub }$ consists of the monomorphisms with cokernel in $\cC_{ \psub } = {}^{ \perp }\cE$, that is, the monomorphisms with semiprojective cokernel; see \ref{bfhpg:Frobenius_E}.
  
  \item  $\fib_{ \psub }$ consists of all epimorphisms.

\end{itemize}

Applying Theorem \ref{thm:Gillespie} to this model category structure gives the following, where the two first items recover the third bullet from \ref{bfhpg:Frobenius_E}, and the last item recovers the first equivalence of categories from \ref{bfhpg:DQA}.
\begin{itemize}
\setlength\itemsep{4pt}
\setcounter{enumi}{3}

  \item  $\cC_{ \psub } \cap \cF_{ \psub } = {}^{ \perp }\cE$ is a Frobenius category.
  
  \item  The class of projective-injective objects is $\cC_{ \psub } \cap \cW_{ \psub } \cap \cF_{ \psub } = {}^{ \perp }\cE \cap \cE = {}_{ Q,A }\!\Prj$, where we used \ref{bfhpg:cotorsion_pair1}.
  
  \item  There is an equivalence of categories 
\[
  \frac{ {}^{ \perp }\cE }{ {}_{ Q,A }\!\Prj }
  \cong \cD_Q( A ),
\]
where the right hand side is defined as
\[
  \cD_Q( A ) =
  \Ho( {}_{ Q,A }\!\Mod ) =
  \weq_{ \psub }^{ -1 }\!{}_{ Q,A }\!\Mod =
  \weq^{ -1 }\!{}_{ Q,A }\!\Mod,
\]  
the homotopy category of the projective model category structure on ${}_{ Q,A }\!\Mod$.
\end{itemize}
\end{bfhpg}

\begin{bfhpg}[The injective model category structure on ${}_{ Q,A }\!\Mod$]
\label{bfhpg:injective_model_category_structure}
By Theorem \ref{thm:Hovey}, the hereditary Hovey triple $( \cC_{ \isub },\cW_{ \isub },\cF_{ \isub } ) = ( {}_{ Q,A }\!\Mod,\cE,\cE^{ \perp } )$ from \ref{bfhpg:Hovey_triple_I} gives a model category structure $( \weq_{ \isub },\cof_{ \isub },\fib_{ \isub } )$ on ${}_{ Q,A }\!\Mod$, called the {\em injective model category structure}, which can be described as follows.
\begin{itemize}
\setlength\itemsep{4pt}

  \item  $\weq_{ \isub }$ consists of the compositions $\pi\iota$ where $\pi$ is an epimorphism with kernel in $\cW_{ \isub } = \cE$ and $\iota$ is a monomorphism with cokernel in $\cW_{ \isub } = \cE$.  Note that
\[  
  \weq_{ \isub } = \weq
\]
by \cite[thm.\ 7.2]{HJ-JLMS}, where $\weq$ is the class from \ref{bfhpg:weq}.

  \item  $\cof_{ \isub }$ consists of all monomorphisms.
    
  \item  $\fib_{ \isub }$ consists of the epimorphisms with kernel in $\cF_{ \isub } = \cE^{ \perp }$, that is, the epimorphisms with semiinjective kernel; see \ref{bfhpg:Frobenius_E}.

\end{itemize}

Applying Theorem \ref{thm:Gillespie} to this model category structure gives the following, where the two first items recover the fourth bullet from \ref{bfhpg:Frobenius_E}, and the last item recovers the second equivalence of categories from \ref{bfhpg:DQA}.
\begin{itemize}
\setlength\itemsep{4pt}
\setcounter{enumi}{3}

  \item  $\cC_{ \isub } \cap \cF_{ \isub } = \cE^{ \perp }$ is a Frobenius category.
  
  \item  The class of projective-injective objects is $\cC_{ \isub } \cap \cW_{ \isub } \cap \cF_{ \isub } = \cE \cap \cE^{ \perp } = {}_{ Q,A }\!\Inj$, where we used \ref{bfhpg:cotorsion_pair3}.
  
  \item  There is an equivalence of categories 
\[
  \frac{ \cE^{ \perp } }{ {}_{ Q,A }\!\Inj }
  \cong
  \cD_Q( A ),
\]
where the right hand side is defined as
\[
  \cD_Q( A ) =
  \Ho( {}_{ Q,A }\!\Mod ) =
  \weq_{ \isub }^{ -1 }\!{}_{ Q,A }\!\Mod =
  \weq^{ -1 }\!{}_{ Q,A }\!\Mod,
\]  
the homotopy category of the injective model category structure on ${}_{ Q,A }\!\Mod$.
\end{itemize}
\end{bfhpg}

\begin{bfhpg}[The projective model category structure on complexes]
\label{bfhpg:projective_model_category_structure_cpx}
Let $Q = \Qcpx$ whence ${}_{ Q,A }\!\Mod = \Ch( A )$; see \ref{bfhpg:examples_of_Q}.  The projective model category structure from \ref{bfhpg:projective_model_category_structure} can be described as follows.
\begin{itemize}
\setlength\itemsep{4pt}

  \item  $\weq_{ \psub } = \weq = \{ \mbox{ quasiisomorphisms } \}$; see the fourth bullet in \ref{bfhpg:E_etc_for_complexes}.
  
  \item  $\cof_{ \psub }$ consists of the monomorphisms whose cokernel is a semiprojective complex; see the second bullet in \ref{bfhpg:projective_model_category_structure} and the first bullet in \ref{bfhpg:Eperp_for_complexes}.
  
  \item  $\fib_{ \psub }$ consists of all epimorphisms; see the third bullet in \ref{bfhpg:projective_model_category_structure}.

\end{itemize}
This is the ``standard model category structure'' on chain complexes described in \cite[def.\ 2.3.3]{Hovey-book}; see \cite[prop.\ 2.3.4, prop.\ 2.3.9]{Hovey-book}.

The corresponding homotopy category is the $Q$-shaped derived category, which for $Q = \Qcpx$ is equivalent to the classic derived category:
\[
  \Ho( {}_{ Q,A }\Mod ) = 
  \cD_Q( A ) \cong 
  \cD( A ),
\]
see \ref{bfhpg:DQA} and \ref{bfhpg:DQA_cpx}.
\end{bfhpg}

\begin{bfhpg}[The injective model category structure on complexes]
\label{bfhpg:injective_model_category_structure_cpx}
Let $Q = \Qcpx$ whence ${}_{ Q,A }\!\Mod = \Ch( A )$.  The injective model category structure from \ref{bfhpg:injective_model_category_structure} can be described as follows.
\begin{itemize}
\setlength\itemsep{4pt}

  \item  $\weq_{ \isub } = \weq = \{ \mbox{ quasiisomorphisms } \}$; see the fourth bullet in \ref{bfhpg:E_etc_for_complexes}.
  
  \item  $\cof_{ \isub }$ consists of all monomorphisms; see the second bullet in \ref{bfhpg:injective_model_category_structure}.

  \item  $\fib_{ \isub }$ consists of the epimorphisms whose kernel is a semiinjective complex; see the third bullet in \ref{bfhpg:injective_model_category_structure} and the second bullet in \ref{bfhpg:Eperp_for_complexes}.

\end{itemize}
This is the ``injective model category structure'' on chain complexes described in \cite[thm.\ 2.3.13]{Hovey-book}.

The corresponding homotopy category is the $Q$-shaped derived category, which for $Q = \Qcpx$ is equivalent to the classic derived category:
\[
  \Ho( {}_{ Q,A }\Mod ) = 
  \cD_Q( A ) \cong 
  \cD( A ),
\]
see \ref{bfhpg:DQA} and \ref{bfhpg:DQA_cpx}.
\end{bfhpg}

\begin{bfhpg}[The projective and injective model category structures on $N$-complexes]
Let $Q = \QNcpx$ whence ${}_{ Q,A }\!\Mod = \ChN( A )$; see \ref{bfhpg:examples_of_Q}.  The projective and injective model category structures from \ref{bfhpg:projective_model_category_structure} and \ref{bfhpg:injective_model_category_structure} can be described analogously to \ref{bfhpg:projective_model_category_structure_cpx} and \ref{bfhpg:injective_model_category_structure_cpx}.  In either case we have
\begin{itemize}
\setlength\itemsep{4pt}

  \item  $\weq = \{ \mbox{ $N$-quasiisomorphisms } \}$; see the last bullet in \ref{bfhpg:E_etc_for_N-complexes}.

\end{itemize}
The corresponding homotopy category is the $Q$-shaped derived category, which for $Q = \QNcpx$ is equivalent to the derived category of $N$-complexes:
\[
  \Ho( {}_{ Q,A }\Mod ) = 
  \cD_Q( A ) \cong 
  \cD_N( A ),
\]
see \ref{bfhpg:DQA} and \ref{bfhpg:DQA_Ncpx}.
\end{bfhpg}

\part{Compact, perfect, and strictly perfect objects in $\cD_Q( A )$}
\label{part:compact}

This part presents some classes of objects and some properties of $\cD_Q( A )$.

\section{Compact, perfect, and strictly perfect objects}
\label{sec:compact}

\begin{bfhpg}[Compact, perfect, and strictly perfect objects of ${}_{ Q,A }\!\Mod$]
\label{bfhpg:perfect}
The full subcategory of {\em compact} objects in $\cD_Q( A )$ is
\begin{itemize}
\setlength\itemsep{4pt}

  \item  
$
  \cD_Q^{ \c }( A ) =
  \{
    C \in \cD_Q( A ) 
    \mid
    \mbox{$\Hom_{ \cD_Q( A ) }( C,- )$ respects set indexed coproducts}
  \},
$
\end{itemize}
see \cite[def.\ 1.6]{Neeman-duality}.

Inspired by \cite[def.\ I.2.1]{SGA6}, the full subcategories of {\em strictly perfect} and {\em perfect} objects in $\cD_Q( A )$ were defined as follows in \cite[def.\ 5.3]{HJ-arXiv}.  The definitions should be read with an understanding that the category $\cD_Q( A ) = \weq^{ -1 }\!{}_{ Q,A }\!\Mod$ has the same objects as ${}_{ Q,A }\!\Mod$.
\begin{itemize}
\setlength\itemsep{6pt}

  \item  
$
  \cD_Q^{ \sperf }( A ) =
  \Bigg\{
    K \in \cD_Q( A ) 
    \:\Bigg|\!
    \begin{array}{l}
      \mbox{the set $\{ q \in Q_0 \mid K( q ) \neq 0 \}$ is finite, and each} \\[1mm] 
      \mbox{$K( q )$ is a finitely generated projective $A$-module} 
    \end{array}
  \Bigg\}
$

  \item  
$
  \cD_Q^{ \perf }( A ) =
  \{
    X \in \cD_Q( A ) 
    \mid
    \mbox{$X \cong K$ in $\cD_Q( A )$ for an object $K \in \cD_Q^{ \sperf }( A ) \cap {}^{ \perp }\cE$}
  \}
$

\end{itemize}
By definition, $\cD_Q^{ \c }( A )$ and $\cD_Q^{ \perf }( A )$ are closed under isomorphisms in $\cD_Q( A )$.

In general, $\cD_Q^{ \sperf }( A )$ is not closed under isomorphisms in $\cD_Q( A )$, and $\cD_Q^{ \perf }( A )$ is not the isomorphism closure of $\cD_Q^{ \sperf }( A )$ in $\cD_Q( A )$ because of the condition $K \in \cD_Q^{ \sperf }( A ) \cap {}^{ \perp }\cE$ in the last bullet above.  However, see Theorem \ref{thm:HJ-arXiv_B}.
\end{bfhpg}

The rationale for the definition of perfect objects is the following two theorems.

\begin{Theorem}[{\cite[thm.\ A]{HJ-arXiv}}]
\label{thm:HJ-arXiv_A}
In general we have $\cD_Q^{ \sperf }( A ) \not\subseteq \cD_Q^{ \c }( A )$.
\end{Theorem}

\begin{Theorem}[{\cite[thm.\ C]{HJ-arXiv}}]
\label{thm:HJ-arXiv_C}
There is an inclusion $\cD_Q^{ \perf }( A ) \subseteq \cD_Q^{ \c }( A )$, which is an equality if and only if $\cD_Q^{ \perf }( A )$ is thick.
\end{Theorem}

Theorem \ref{thm:HJ-arXiv_C} motivates the following conjecture, which is well known to be true for $\cD( A )$.

\begin{Conjecture}
There is an equality $\cD_Q^{ \perf }( A ) = \cD_Q^{ \c }( A )$.
\end{Conjecture}

In some cases, the relation between the categories $\cD_Q^{ \c }( A )$, $\cD_Q^{ \sperf }( A )$, and $\cD_Q^{ \perf }( A )$ simplifies as described by the following theorem.  It applies in particular to $\cD( A )$, which can be obtained as $\cD_Q( A )$ for the category $Q = \Qcpx$ which has no cycles.

\begin{Theorem}[{\cite[thm.\ B]{HJ-arXiv}}]
\label{thm:HJ-arXiv_B}
Assume that $Q$ has no cycles in the sense of \ref{bfhpg:cycle} or that the left global dimension of $A$ is finite.  Then the following hold.
\begin{itemize}
\setlength\itemsep{4pt}

  \item  $\cD_Q^{ \sperf }( A ) \subseteq \cD_Q^{ \c }( A )$.

  \item  $\cD_Q^{ \perf }( A )$ is the isomorphism closure of $\cD_Q^{ \sperf }( A )$ in $\cD_Q( A )$. 

\end{itemize}
\end{Theorem}

Finally, $\cD_Q( A )$ always enjoys the following good property.

\begin{Theorem}[{\cite[thm.\ D]{HJ-arXiv}}]
\label{thm:HJ-arXiv_D}
The category $\cD_Q( A )$ is compactly generated in the sense of \cite[def.\ 1.7]{Neeman-duality}.  A set of compact generators is given by $\{ S_q( A ) \mid q \in Q_0 \}$, where the functor $S_q( A ) \in {}_{ Q,A }\!\Mod$ is defined by
\[
  \big( S_q( A ) \big)( - ) = \big( S \langle q \rangle ( - ) \big) \Tensor{ \Bk } A,
\]  
see \ref{bfhpg:the_functors_H}.
\end{Theorem}

If $Q$ is $\Qcpx$ or $\QNcpx$, then ${}_{ Q,A }\!\Mod$ is $\Ch( A )$ or $\ChN( A )$, and $S_q( A )$ is the complex or $N$-complex which has $A$ placed in degree $q$ and zeroes elsewhere.  If $Q = \Qcpx$ then $\cD_Q( A ) = \cD( A )$, see \ref{bfhpg:DQA_cpx}, and this set of compact generators is well known.

\appendix

\part{Appendices}
\label{part:compendium}
\renewcommand{\thesection}{\Roman{section}}

This part contains two appendices on some key classes of categories: Frobenius, triangulated, and abelian model categories.

\section{Frobenius and triangulated categories} 
\label{app:I}

\begin{bfhpg}[Exact categories]
\label{bfhpg:exact}
The notion of {\em exact category} was introduced by Quillen; see \cite[sec.\ 2]{Quillen-K} or \cite[def.\ 2.1]{Buehler} for a more recent exposition.  We will not reproduce the definition in full but merely say the following.
\begin{itemize}
\setlength\itemsep{4pt}

  \item  An exact category is pair $( \cF,\cS )$ where $\cF$ is an additive category, $\cS$ a class of so-called {\em conflations}.
  
  \item  Each conflation is a diagram of the form $f' \xrightarrow{ \varphi' } f \xrightarrow{ \varphi } f''$ where $\varphi'$ is a kernel of $\varphi$ and $\varphi$ is a cokernel of $\varphi'$.
  
  \item  The conflations are subject to a list of axioms.

\end{itemize}
The canonical example is that $\cF$ is an extension closed subcategory of an abelian category $\cA$, and that $\cS$ is the class of short exact sequence in $\cA$ which have each term in $\cF$.  In this case $\cS$ is implicit in $\cF$, and we abuse terminology by saying that ``$\cF$ is an exact category''.
\end{bfhpg}

\begin{bfhpg}[Frobenius categories]
\label{bfhpg:Frobenius}
Let $( \cF,\cS )$ be an exact category.
\begin{itemize}
\setlength\itemsep{4pt}

  \item  An object $p \in \cF$ is {\em projective} if $\cF( p,- )$ maps conflations to short exact sequences.
  
  \item  There are {\em enough projective objects} if each $f \in \cF$ permits a conflation $f' \xrightarrow{} p \xrightarrow{} f$ with $p$ projective.
  
  \item  There are dual definitions of {\em injective} objects and of {\em enough injective objects}.
  
  \item  We say that $( \cF,\cS )$ is a {\em Frobenius category} if it has enough projective and enough injective objects and the projective and injective objects coincide; see \cite[sec.\ 2.1]{Happel-book}.  These are then referred to as {\em projective-injective} objects. 

\end{itemize}
As above, if $\cF$ is given as an extension closed subcategory of an abelian category $\cA$, and $\cS$ is the class of short exact sequence in $\cA$ which have each term in $\cF$, then we abuse terminology by saying that ``$\cF$ is a Frobenius category''.

A simple example is that $\cF$ is $\mod( \Lambda )$, the category of finitely generated $\Lambda$-left modules where $\Lambda$ is a finite dimensional self injective algebra over a field $k$, and that $\cS$ is the class of all short exact sequences in $\cF$.
\end{bfhpg}

\begin{bfhpg}[Triangulated categories]
\label{bfhpg:triangulated}
An early reference for the definition of triangulated categories is \cite[sec.\ I.1]{Hartshorne-resdual}; see also \cite[sec. I.1.1]{Happel-book}.  We will not reproduce the definition in full but merely say the following.
\begin{itemize}
\setlength\itemsep{4pt}

  \item  A triangulated category is a triple $( \cD,\Sigma,\Delta )$ where $\cD$ is an additive category, $\Sigma$ an automorphism of $\cD$ called the {\em suspension functor}, and $\Delta$ a class of so-called {\em triangles}.
  
  \item  Each triangle is a diagram of the form $d' \xrightarrow{ \delta' } d \xrightarrow{ \delta } d'' \xrightarrow{} \Sigma d'$ where $\delta'$ is a weak kernel of $\delta$ and $\delta$ is a weak cokernel of $\delta'$.
  
  \item  The triangles are subject to a list of axioms.

\end{itemize}
A main example is $\cD( A )$; see \cite[sec.\ I.4]{Hartshorne-resdual} and \cite[sec.\ I.3.3]{Happel-book}.
\end{bfhpg}

\begin{bfhpg}[The stable category of a Frobenius category]
\label{bfhpg:stable}
Let $( \cF,\cS )$ be a Frobenius category and let $\cP$ be the full subcategory of projective-injective objects.  The {\em stable} category $\frac{ \cF }{ \cP }$ is the naive quotient with the same objects as $\cF$ and morphisms
\[
  \Hom_{ \frac{ \cF }{ \cP } }( f',f )
  =
  \frac{ \Hom_{ \cF }( f',f ) }{ \{ \mbox{ morphisms which factorise through an object of $\cP$ } \} }.
\]

The stable category is a triangulated category in a canonical way by \cite[thm.\ I.2.6]{Happel-book}.  We will not provide full details but merely say the following.
\begin{itemize}
\setlength\itemsep{4pt}

  \item  A conflation $f' \xrightarrow{} f \xrightarrow{} f''$ in $\cF$ induces a triangle $f' \xrightarrow{} f \xrightarrow{} f'' \xrightarrow{} \Sigma f'$ in $\frac{ \cF }{ \cP }$.

  \item  Up to isomorphism, the suspension $\Sigma f$ of an object $f$ is $\Sigma f = f'$ where $f \xrightarrow{} p \xrightarrow{} f'$ is a conflation with $p$ projective-injective.

\end{itemize} 
\end{bfhpg}

\section{Abelian model categories}
\label{app:II}

\begin{bfhpg}[Model categories]
\label{bfhpg:model_categories}
The notion of {\em model category} was introduced by Quillen; see \cite[sec.\ I.1]{Quillen-hot}.  In this paper, we will use the updated definition of \cite[def.\ 1.1.3]{Hovey-book}.  We will not reproduce the definition in full but merely say the following.
\begin{itemize}
\setlength\itemsep{4pt}

  \item  A model category is a quadruple $(\cA,\weq,\cof,\fib)$, where $\cA$ is a category and $(\weq,\cof,\fib)$ are three classes of morphisms called the {\em weak equivalences}, {\em cofibrations}, and {\em fibrations}.  The triple $( \weq,\cof,\fib )$ is often referred to as a {\em model category structure} on $\cA$.  
  
  \item  A {\em trivial} (co)fibration is a weak equivalence which is also a (co)fibration.
  
  \item  Using $( \weq,\cof,\fib )$, one can define the so-called {\em cofibrant}, {\em fibrant}, {\em trivally cofibrant}, and {\em trivially fibrant} objects.  Objects which are both cofibrant and fibrant are called {\em cofibrant-fibrant}; they form the full subcategory $\cA_{ \cf }$.  
  
  \item  There is a notion of homotopy, denoted $\sim$, which is an equivalence relation on each $\Hom$ set between cofibrant-fibrant objects, compatible with composition of morphisms.

\end{itemize}
The {\em homotopy category} is
\[
  \Ho( \cA ) = \weq^{ -1 }\!\cA,
\]
obtained from $\cA$ by formally inverting each weak equivalence.  Quillen's fundamental theorem of model categories, \cite[thm.\ 1']{Quillen-hot}, states that the inclusion $\cA_{ \cf } \hookrightarrow \cA$ induces an equivalence of categories
\[
  \cA_{ \cf }/\!\sim \;\: \cong \: \Ho( \cA ).
\]
There are at least two points to this: First, it shows that $\weq^{ -1 }\!\cA$ exists without set theoretical issues, which could otherwise encumber its construction.  Secondly, concrete computations may be more feasible in $\cA_{ \cf }/\!\sim$ than in $\weq^{ -1 }\!\cA$.  
\end{bfhpg}

\begin{bfhpg}[Abelian model categories]
\label{bfhpg:abelian_model_categories}
An {\em abelian model category}, as defined by Hovey \cite[def.\ 2.1]{Hovey-MathZ}, is a model category $(\cA,\weq,\cof,\fib)$ where $\cA$ is an abelian category, such that the following hold.
\begin{itemize}
\setlength\itemsep{4pt}

  \item  A morphism is a (trivial) cofibration if and only if it is a monomorphism with (trivially) cofibrant cokernel.
  
  \item  A morphism is a (trivial) fibration if and only if it is an epimorphism with (trivially) fibrant kernel.

\end{itemize}
\end{bfhpg}

\begin{bfhpg}[Wide subcategories]
\label{bfhpg:wide}
Let $\cA$ be an abelian category, $\cW$ a full subcategory.
\begin{itemize}
\setlength\itemsep{4pt}

  \item  $\cW$ is {\em wide} if it is closed under summands and has the two-out-of-three property; that is, given a short exact sequence, if two of the objects are in the subcategory then so is the third.
  
\end{itemize}
\end{bfhpg}

\begin{bfhpg}[Cotorsion pairs]
\label{bfhpg:cotorsion_pair}
Let $\cA$ be an abelian category, $( \cX,\cY )$ a pair of full subcategories.  Recall that $\perp$ denotes perpendicular full subcategories with respect to $\Ext^1$; see \ref{bfhpg:QAMod}.
\begin{itemize}
\setlength\itemsep{4pt}

  \item  $( \cX,\cY )$ is a {\em cotorsion pair} if $\cX^{ \perp } = \cY$ and $\cX = {}^{ \perp }\cY$.  

  \item  A cotorsion pair $( \cX,\cY )$ is {\em functorially complete} if each $a \in \cA$ permits short exact sequences
\[
  0 \xrightarrow{} y \xrightarrow{} x \xrightarrow{} a \xrightarrow{} 0
  \;\;,\;\;
  0 \xrightarrow{} a \xrightarrow{} y' \xrightarrow{} x' \xrightarrow{} 0,
\]
which depend functorially on $a$ and have $x,x' \in \cX$ and $y,y' \in \cY$.  

  \item  A cotorsion pair $( \cX,\cY )$ is {\em hereditary} if $\Ext_{ \cA }^{\geqslant 1}( \cX,\cY ) = 0$.  

\end{itemize}
The notion of cotorsion pair is due to Salce, \cite[p.\ 12]{Salce}.  See also \cite[def.\ 2.2.1, lem.\ 2.2.6, lem.\ 2.2.10]{Goebel-Trlifaj-book}, \cite[def.\ 2.3]{Hovey-MathZ}.
\end{bfhpg}

\begin{bfhpg}[Hovey triples]
\label{bfhpg:Hovey_triple}
A {\em Hovey triple} in an abelian category is a triple $( \cC,\cW,\cF )$ of full subcategories such that
\begin{itemize}
\setlength\itemsep{4pt}

  \item  $\cW$ is a wide subcategory,

  \item  $( \cC,\cW \cap \cF )$ and $( \cC \cap \cW,\cF )$ are functorially complete cotorsion pairs.

\end{itemize}
A Hovey triple is called {\em hereditary} if the cotorsion pairs $( \cC,\cW \cap \cF )$ and $( \cC \cap \cW,\cF )$ are hereditary.  See \cite[thm.\ 2.2 and sec.\ 4]{Gillespie-BLMS}, \cite[thm.\ 2.2]{Hovey-MathZ}.
\end{bfhpg}

\begin{Theorem}[Hovey; see {\cite[thm.\ 2.2 and def.\ 5.1]{Hovey-MathZ}}]
\label{thm:Hovey}
Let $\cA$ be an abelian category with set indexed limits and colimits and consider:
\begin{enumerate}
\setlength\itemsep{4pt}

  \item  The class of abelian model categories of the form $(\cA,\weq,\cof,\fib)$,
  
  \item  The class of Hovey triples $( \cC,\cW,\cF )$ in $\cA$. 

\end{enumerate}
There is a bijection between (i) and (ii).  If $(\cA,\weq,\cof,\fib)$ and $( \cC,\cW,\cF )$ correspond under the bijection, then on the one hand,
\begin{itemize}
\setlength\itemsep{4pt}

  \item  $\cC$ consists of the cofibrant objects, 

  \item  $\cC \cap \cW$ consists of the trivially cofibrant objects,
  
  \item  $\cF$ consists of the fibrant objects,

  \item  $\cW \cap \cF$ consists of the trivially fibrant objects,
  
  \item  $\cW$ consists of the trivial objects, that is, the objects $w$ such that $0 \xrightarrow{} w$ is a weak equivalence.

\end{itemize}
On the other hand,
\begin{itemize}
\setlength\itemsep{4pt}

  \item  $\weq$ consists of the compositions $\pi\iota$ where $\pi$ is an epimorphism with kernel in $\cW$ and $\iota$ is a monomorphism with cokernel in $\cW$\!,

  \item  $\cof$ consists of the monomorphisms with cokernel in $\cC$,
  
  \item  $\fib$ consists of the epimorphisms with kernel in $\cF$.

\end{itemize}
\end{Theorem}

\begin{Remark}
\label{rmk:Acf}
Note that in Theorem \ref{thm:Hovey} the cofibrant-fibrant objects are
\begin{itemize}
\setlength\itemsep{4pt}

  \item  $\cA_{ \cf } = \cC \cap \cF$.  

\end{itemize}
\end{Remark}

\begin{Theorem}[Gillespie; see {\cite[thm.\ 2.6(i) and prop.\ 4.2]{Gillespie-BLMS}}]
\label{thm:Gillespie}
Let $\cA$ be an abelian category with set indexed limits and colimits, $( \cC, \cW, \cF )$ a hereditary Hovey triple in $\cA$, and $( \cA,\weq,\cof,\fib )$ the corresponding abelian model category under Theorem \ref{thm:Hovey}.  Then
\begin{itemize}
\setlength\itemsep{4pt}

  \item  $\cA_{ \cf } = \cC \cap \cF$ is a Frobenius category.  The conflations are the short exact sequences in $\cA$ which have each term in $\cC \cap \cF$,
  
  \item  The class of projective-injective objects of $\cC \cap \cF$ is $\cC \cap \cW \cap \cF$,
  
  \item  The inclusion $\cA_{ \cf } = \cC \cap \cW \hookrightarrow \cA$ induces an equivalence of categories
\[
  \frac{ \cC \cap \cF }{ \cC \cap \cW \cap \cF }
  \cong \weq^{ -1 }\!\cA = \Ho( \cA ).
\]

\end{itemize}
In particular, $\Ho( \cA )$ is triangulated; see \ref{bfhpg:stable}.  
\end{Theorem}

The final bullet can be viewed as an instance of Quillen's fundamental theorem of model categories; see \ref{bfhpg:model_categories} and \cite[thm.\ 1']{Quillen-hot}.

\medskip
\noindent
{\bf Acknowledgement.}
The scientific and organising committee of the Abel Symposium 2022 consisted of Paul Balmer, Petter Andreas Bergh, Bernhard Keller, Henning Krause, Steffen Oppermann, and \O yvind Solberg.  We thank them for the invitation to give a talk on this material and to submit this paper to the proceedings of the symposium.

Part of this paper was written when the second author was a Fellow at the Centre for Advanced Study of the Norwegian Academy of Science and Letters.  He thanks the Centre for its hospitality and the organisers, Aslak Buan and Steffen Oppermann, for the invitation to attend the programme ``Representation Theory: Combinatorial Aspects and Applications''.

This work was supported by a DNRF Chair from the Danish National Research Foundation (grant DNRF156), by a Research Project 2 from the Independent Research Fund Denmark (grant 1026-00050B), and by Aarhus University Research Foundation (grant AUFF-F-2020-7-16).


\begin{thebibliography}{39}

\bibitem{ASS}  I.\ Assem, D.\ Simson, and A.\ Skowro\'{n}ski, ``Elements of the representation theory of associative algebras'', Vol.\ 1, London Math.\ Soc.\ Stud.\ Texts, Vol.\ 65, Cambridge University Press, Cambridge, 2006.




\bibitem{Bahiraei}  P.\ Bahiraei, {\it Cotorsion pairs and adjoint functors in the homotopy category of $N$-complexes}, J.\ Algebra Appl.\ {\bf 19} (2020), 2050236.

\bibitem{SGA6}  P.\ Berthelot, A.\ Grothendieck, and L.\ Illusie, ``Th\'{e}orie des intersections et th\'{e}or\`{e}me de Riemann--Roch (SGA 6)'', in cooperation with D.\ Ferrand, J.\ P.\ Jouanolou, O.\ Jussila, S.\ Kleiman, M.\ Raynaud, and J.\ P.\ Serre, Lecture Notes in Math., Vol.\ 225, Springer, Berlin--Heidelberg--New York, 1971.


\bibitem{Boekstedt-Neeman}  M.\ B\"{o}kstedt and A.\ Neeman, {\it Homotopy limits in triangulated categories}, Compositio Math.\ {\bf 86} (1993), 209--234. 

\bibitem{Buehler}  T.\ B\"{u}hler, {\it Exact categories}, Expo.\ Math.\ {\bf 28} (2010), 1--69.


\bibitem{DSS}  I.\ Dell'Ambrogio, G.\ Stevenson, and J.\ \v{S}\v{t}ov\'{\i}\v{c}ek, {\it Gorenstein homological algebra and universal coefficient theorems}, Math.\ Z.\ {\bf 287} (2017), 1109--1155.








\bibitem{Garcia-Rozas-book}  J.\ R.\ Garc\'{\i}a Rozas, ``Covers and envelopes in the category of complexes of modules'', Chapman \& Hall/CRC Res.\ Notes Math., Vol. 407, Chapman \& Hall/CRC, Boca Raton, Florida, 1999.

\bibitem{Gillespie-BLMS}  J.\ Gillespie, {\it Hereditary abelian model categories}, Bull.\ London Math.\ Soc.\ {\bf 48} (2016), 895--922.

\bibitem{Gillespie-construct}  J.\ Gillespie, {\it How to construct a Hovey triple from two cotorsion pairs}, Fund.\ Math.\ {\bf 230} (2015), 281--289.

\bibitem{Gillespie-exact}  J.\ Gillespie, {\it Model structures on exact categories}, J.\ Pure Appl.\ Algebra {\bf 215} (2011), 2892--2902.

\bibitem{Gillespie-flat}  J.\ Gillespie, {\it The flat model structure on {${\rm Ch}(R)$}}, Trans.\ Amer.\ Math.\ Soc.\ {\bf 356} (2004), 3369--3390.


\bibitem{Goebel-Trlifaj-book}  R.\ G{\"o}bel and J.\ Trlifaj, ``Approximations and endomorphism algebras of modules'', Vol.\ 41, Walter de Gruyter, Berlin, 2006.

\bibitem{Grothendieck-Tohoku}  A.\ Grothendieck, {\it Sur quelques points d'alg\`{e}bre homologique}, Tohoku Math.\ J.\ (2) {\bf 9} (1957) 119--221.

\bibitem{Happel-book}  D.\ Happel, ``Triangulated categories in the representation theory of finite dimensional algebras'', London Math.\ Soc.\ Lecture Note Ser., Vol.\ 119, Cambridge University Press, Cambridge, 1988.

\bibitem{Hartshorne-resdual}  R.\ Hartshorne, ``Residues and duality'', lecture notes of a seminar on the work of A.\ Grothendieck, given at Harvard 1963/64.  With an appendix by P.\ Deligne, Lecture Notes in Math., Vol.\ 20, Springer, Berlin--New York, 1966.





\bibitem{HJ-Model_cats}  H.\ Holm and P.\ J\o rgensen, {\it Model categories of quiver representations}, Adv.\ Math.\ {\bf 357} (2019), 106826. 

\bibitem{HJ-JLMS}  H.\ Holm and P.\ J\o rgensen, {\it The $Q$-shaped derived category of a ring}, J.\ London Math.\ Soc.\ (2) {\bf 106} (2022), 3263--3316.

\bibitem{HJ-arXiv}  H.\ Holm and P.\ J\o rgensen, The $Q$-shaped derived category of a ring -- compact and perfect objects, preprint (2022).  {\tt arXiv:2208.13282. }

\bibitem{Hovey-MathZ}  M.\ Hovey, {\it Cotorsion pairs, model category structures, and representation theory}, Math.\ Z.\ {\bf 241} (2002), 553--592.

\bibitem{Hovey-book}  M.\ Hovey, ``Model categories'', Math.\ Surveys Monogr., Vol.\ 63, American Mathematical Society, Providence, 1999.

\bibitem{Iyama-Kato-Miyachi-N}  O.\ Iyama, K.\ Kato, and J. Miyachi, {\it Derived categories of $N$-complexes}, J.\ London Math.\ Soc.\ (2) {\bf 96} (2017), 687--716.

\bibitem{Iyama-Minamoto-1}  O.\ Iyama and H.\ Minamoto, $\cA$-derived categories, in preparation.

\bibitem{Iyama-Minamoto-2}  O.\ Iyama and H.\ Minamoto, {\it On a generalization of complexes and their derived categories}, extended abstract of a talk at the 47th Ring and Representation Theory Symposium, Osaka City University, 2014.


\bibitem{Kapranov}  M.\ M.\ Kapranov, On the $q$-analog of homological algebra, preprint (1996).  {\tt arXiv:9611005v1. }









\bibitem{Neeman-duality}  A.\ Neeman, {\it The Grothendieck duality theorem via Bousfield's techniques and Brown representability}, J.\ Amer.\ Math.\ Soc.\ {\bf 9} (1996), 205--236.


\bibitem{Oberst-Roehrl}  U.\ Oberst and H.\ Rohrl, {\it Flat and coherent functors}, J.\ Algebra {\bf 14} (1970), 91--105.



\bibitem{Quillen-K}  D.\ Quillen, {\it Higher algebraic $K$-theory: I}, pp.\ 85--147 in ``Algebraic $K$-theory, I: Higher $K$-theories'' (proceedings of the conference held at the Seattle Research Center of the Battelle Memorial Institute, 28 August to 8 September, 1972, edited by H.\ Bass), Lecture Notes in Math., Vol.\ 341, Springer, Berlin, 1973.

\bibitem{Quillen-hot}  D.\ Quillen, ``Homotopical algebra'', Lecture Notes in Math., Vol.\ 43, Springer, Berlin--Heidelberg, 1967.


\bibitem{Salce}  L.\ Salce, {\it Cotorsion theories for abelian groups}, pp.\ 11-32 in ``Symposia Mathematica, Vol.\ {XXIII} (Conf.\ Abelian Groups and their Relationship to the Theory of Modules, INDAM, Rome, 1977)'', Academic Press, London--New York, 1979.



\bibitem{Schiffler}  R.\ Schiffler, ``Quiver representations'', CMS Books Math., Springer, Cham, 2014.


\bibitem{Stenstroem}  B.\ Stenstr\"{o}m, ``Rings of quotients'', Grundlehren Math.\ Wiss., Vol.\ 217, Springer, Berlin--Heidelberg--New York, 1975.

\bibitem{Stovicek-exact}  J.\ \v{S}\v{t}ov\'{\i}\v{c}ek, {\it Exact model categories, approximation theory, and cohomology of quasi-coherent sheaves}, pp.\ 297-367 in ``Advances in representation theory of algebras'' (proceedings of the International Conference on Representations of Algebras, Bielefeld 2015, edited by Benson, Krause, and Skowro\'{n}ski), EMS Series of Congress Reports, European Mathematical Society, Z\"{u}rich, 2013.


\end{thebibliography}
\end{document}